\documentclass[pdflatex,sn-mathphys-num]{sn-jnl}

\usepackage{amsmath,amssymb,amsthm,mathtools,amsfonts}
\usepackage{graphicx}
\usepackage{multirow}
\usepackage[title]{appendix}
\usepackage{xcolor}
\usepackage{textcomp}
\usepackage{manyfoot}
\usepackage{booktabs}
\usepackage{lineno}
\usepackage{tabularx}
\usepackage[ruled,vlined,algo2e]{algorithm2e}
\usepackage{listings}
\usepackage{hyperref}
\usepackage{enumitem}

\hypersetup{
  colorlinks=true,
  linkcolor=blue!55!black,
  citecolor=blue!55!black,
  urlcolor=blue!55!black
}

\usepackage{enumitem}
\usepackage{float}

\newtheorem{theorem}{Theorem}[section]
\newtheorem{proposition}[theorem]{Proposition}
\newtheorem{lemma}[theorem]{Lemma}

\theoremstyle{definition}

\newcommand{\R}{\mathbb{R}}
\newcommand{\Jell}{\mathcal{J}_{\ell}}
\newcommand{\Ffun}{\mathcal{F}}
\newcommand{\eps}{\varepsilon}
\newcommand{\dd}{\,\mathrm{d}}
\newcommand{\norm}[1]{\left\lVert #1\right\rVert}
\newcommand{\abs}[1]{\left\lvert #1\right\rvert}

\numberwithin{equation}{section}

\begin{document}

\title{A Nonlocal Perfusion-Gated Angiogenesis System: Global Weak Solutions, Fast-Signal Limits, and Invasion Fronts}

\author[1]{\fnm{Jiguang} \sur{Yu}}\email{jyu678@bu.edu}
\equalcont{These authors contributed equally to this work as co-first authors.}

\author*[2]{\fnm{Louis Shuo} \sur{Wang}}\email{wang.s41@northeastern.edu}
\equalcont{These authors contributed equally to this work as co-first authors.}

\affil[1]{\orgdiv{College of Engineering},
  \orgname{Boston University},
  \orgaddress{\city{Boston}, \postcode{02215}, \state{MA}, \country{USA}}}

\affil[2]{\orgdiv{Department of Mathematics},
  \orgname{Northeastern University},
  \orgaddress{\city{Boston}, \postcode{02115}, \state{MA}, \country{USA}}}

\abstract{
Continuum angiogenesis systems often use local endothelial or vessel density
as an oxygen-delivery proxy without distinguishing structurally formed vessels from pressure-supported vascular function.  We formulate a two-dimensional model
in which lumenized density determines a normalized nonlocal conductivity, a
globally solved pressure field, and a flow-functional density that gates
oxygen delivery and vessel regression.  For fixed positive regularization
parameters and signal relaxation times, we prove Lipschitz stability of the
pressure--perfusion map
and global bounded weak solvability; strong vessel compactness follows from a
nonlocal ordinary differential equation stability estimate rather than spatial smoothing.  As the oxygen
and vascular endothelial growth factor timescales vanish simultaneously, weak solutions converge along a
subsequence to a parabolic--elliptic--ordinary differential equation system, with both fast fields
converging strongly in \(L^2(0,T;H^1)\).  A locally frozen one-dimensional
reduction admits monotone fronts precisely at or above the threshold
\(B+2\sqrt{DR}\) for \(B\geq-\sqrt{DR}\) and \(-DR/B\) otherwise.
Finite-volume experiments verify the analytical mechanisms and show that
equal-mass vessel fields can generate distinct functional masses and
oxygenation levels.}

\keywords{Angiogenesis; nonlocal perfusion; chemotaxis--haptotaxis;
global weak solutions; fast-signal limit; traveling waves.}

\pacs[MSC classification]{35K55, 35D30, 35B25, 35C07, 35Q92, 92C17}

\maketitle

\section{Introduction}
\label{sec:introduction}
\begin{figure}[htbp]
\centering
\includegraphics[width=\linewidth]{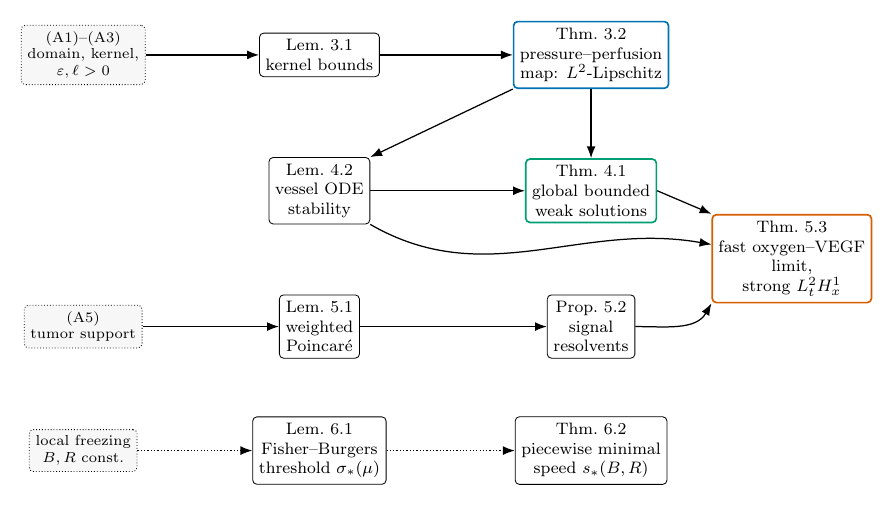}
\caption{Logical architecture of the results.  Solid arrows indicate ``is
used in the proof of''; dotted boxes are assumption blocks; the dotted
track is the locally frozen scalar reduction of Section~\ref{sec:waves},
which adds no global partial differential equation assumption.  The mechanisms along the arrows:
Theorem~\ref{thm:perfusion} supplies the Lipschitz elliptic closure;
Lemma~\ref{lem:vascular} replaces the unavailable spatial compactness of
the vessel field in Theorems~\ref{thm:global} and~\ref{thm:fast};
Lemma~\ref{lem:weighted-poincare} gives the oxygen coercivity behind the
Minty and energy-pairing identification in Theorem~\ref{thm:fast};
Lemma~\ref{lem:burgers} is the variational threshold behind
Theorem~\ref{thm:wave}.}
\label{fig:roadmap}
\end{figure}

\begin{figure}[htbp]
\centering
\includegraphics[width=\linewidth]{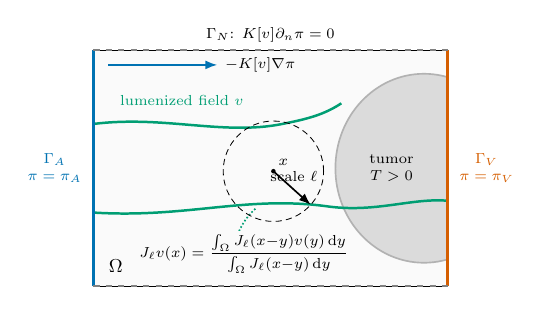}
\caption{Domain geometry.  Pressure is imposed on the arterial and venous
parts \(\Gamma_A,\Gamma_V\) and insulated on \(\Gamma_N\); the tumor
profile \(T\) is prescribed.  The lumenized field \(v\) enters the
conductivity only through the boundary-normalized average \(\Jell v\)
over a kernel neighborhood at scale \(\ell\), whose normalization by
\(\int_\Omega J_\ell(x-y)\dd y\) keeps the average well defined and
order-preserving up to the boundary (Lemma~\ref{lem:kernel}).}
\label{fig:domain}
\end{figure}

Continuum tumor angiogenesis models commonly describe endothelial migration by diffusion, chemotaxis toward a tumor-derived signal, and haptotaxis along the extracellular matrix
\cite{anderson1998continuous,hillen2009user,vilanova2017computational,
de2025numerical}.
Oxygen-regulated extensions add an important metabolic feedback, but may use
local endothelial density as the vascular oxygen source
\cite{de2025analytical}.  Conversely, hybrid discrete and
mixed-dimensional models calculate intravascular pressure, tissue exchange,
and functional perfusion on explicit networks
\cite{mcdougall2006mathematical,berrone2023optimization,duswald2024bridging,
hadjigeorgiou2024hybrid}, building on flow-adaptation and
network-oxygen-transport principles established for microvasculature
\cite{pries1998structural,secomb2004green}.
We study a continuum alternative in which motile endothelial occupancy, formed vascular structural density, and flow-functional density are separate fields. The lumenized field is averaged by a boundary-normalized kernel
--- a nonlocal construction familiar from continuum cell-adhesion models
\cite{armstrong2006continuum}, used here inside an elliptic coefficient rather
than a flux --- and inserted into a uniformly elliptic conductivity.  A macroscopic
arterial--venous pressure solve produces a flow proxy, and the functional
density \(\Ffun[v]=vG(q[v])\) weights both oxygen delivery and low-flow
regression.  This sacrifices exact graph topology in exchange for a closure
that can be analyzed inside the evolutionary partial differential equation system.  Related elliptic
coefficient maps and global/local perfusion couplings provide analytical
context \cite{griepentrog2001linear,felisi2024full,bociu2026existence}, but they
do not couple an evolving lumenized conductivity to angiogenesis dynamics in
the manner studied here.
The proof has two nonstandard elements.  First, the nondiffusive vessel field
does not gain compactness from the Aubin--Lions lemma, and the direct factor
\(v\) prevents a compactness claim for \(\Ffun\).  We instead combine
Lipschitz stability of the elliptic closure with a quantitative vessel-ordinary differential equation (ODE) estimate.  Second, in the
simultaneous fast oxygen--vascular endothelial growth factor (VEGF) limit, neither fast time derivative is
uniformly controlled.  Compactness is first obtained for the slow variables;
the fast fields are then identified and shown to converge strongly using
monotone elliptic resolvents and dynamic energy pairings.  These issues distinguish the
analysis from local Keller--Segel and chemotaxis--haptotaxis theory
\cite{winkler2010aggregation,tao2014energy,bellomo2015toward,chen2024negligibility,
rani2024quasilinear,huo2024global} and from prescribed-source
fast-signal limits
\cite{mizukami2019fast,wang2019fast,li2021convergence,mizukami2018fast}.
The paper also derives a local scalar invasion diagnostic. Exact
reaction--diffusion--convection theory
\cite{malaguti2002travelling,benguria2004minimal,drabek2026traveling}
(see \cite{wang2013mathematics} for a survey of traveling waves in chemotaxis
systems) shows that the usual leading-edge formula fails to give the exact
minimal speed under sufficiently strong adverse taxis.  Computations are used
only to verify the analytical limits and isolate the effect of flow weighting;
they do not constitute parameter calibration or clinical prediction.
The principal contributions are: (i) well-posedness, bounds, measurability,
and \(L^2\)-Lipschitz stability of the pressure--flow--functional map;
(ii) global bounded weak solutions obtained with a vessel-stability closure;
(iii) a fully coupled simultaneous oxygen--VEGF fast-signal limit with strong
\(L^2_tH^1_x\) convergence; and (iv) a complete piecewise minimal-speed
theorem for the locally frozen Fisher--Burgers reduction.  The parameters
\(\eps>0\) and \(\ell>0\) remain fixed, the tumor is prescribed.  
Figure~\ref{fig:roadmap} summarizes the logical architecture of the results and the mechanism each step contributes. The geometric setting --- the mixed arterial--venous--insulated boundary decomposition and the boundary-normalized kernel neighborhood --- is sketched in Figure~\ref{fig:domain}.

The remainder of this paper is organized as follows. Section~\ref{sec:model} details the model derivation, biological balances, nondimensional scaling, and standing assumptions. Section~\ref{sec:perfusion} establishes the well-posedness and Lipschitz stability of the nonlocal pressure--perfusion map. In Section~\ref{sec:global}, we prove the existence of global bounded weak solutions, utilizing a vascular stability closure. Section~\ref{sec:fast} analyzes the simultaneous fast oxygen--VEGF limit and demonstrates the strong convergence properties of the fast fields. Section~\ref{sec:waves} investigates a locally frozen Fisher--Burgers reduction to derive exact minimal invasion-front speeds. Section~\ref{sec:numerical} presents numerical simulations that verify our analytical results and highlight the effects of flow weighting. Finally, Section~\ref{sec:discussion} concludes with a discussion of the model's limitations and future directions. 

\section{Model derivation and assumptions}
\label{sec:model}
\subsection{Dimensional balances and biological levels}
\begin{figure}[htbp]
\centering
\includegraphics[width=\linewidth]{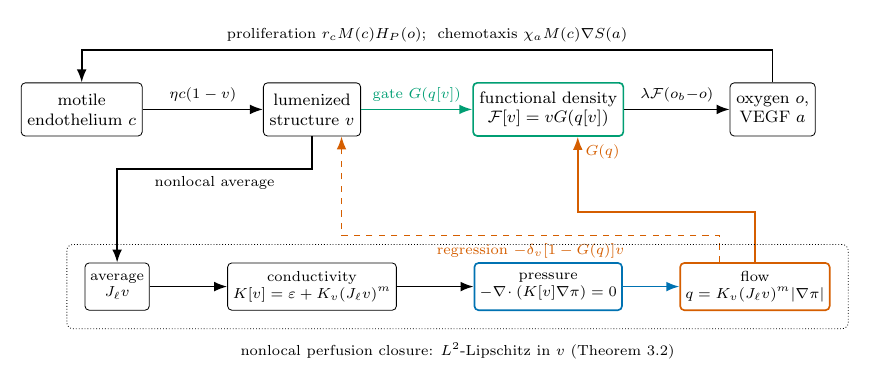}
\caption{Model architecture.  Top row: the biological hierarchy
\eqref{eq:levels} from motile endothelium through lumenized structure to
the flow-functional density that gates oxygen delivery.  Bottom row: the
nonlocal perfusion closure --- kernel average, uniformly elliptic
conductivity, global pressure solve, flow magnitude --- whose
\(L^2\)-Lipschitz continuity (Theorem~\ref{thm:perfusion}) is the
interface used by every compactness argument.  The dashed arrow is the
low-flow regression feedback in \eqref{eq:v}.}
\label{fig:architecture}
\end{figure}

Let \(X\in\widehat\Omega\subset\R^2\), and let \(\mathsf t\) denote
dimensional time.
The variables \(C,P,F,A,V,O,\Pi\) denote, respectively, endothelial areal
density, effective net proteolytic activity, matrix areal density, VEGF areal
concentration, vascular structure length per unit tissue area, oxygen areal
concentration, and vascular pressure. 
Vessel segments are interpreted as lumenized and the vascular structural density $V$ represents the continuum density of vessel segments with an established lumen, but not necessarily effective perfusion.
Their units are
\([C]=N_cL^{-2}\), \([P]=N_pL^{-2}\), \([F]=M_fL^{-2}\),
\([A]=N_aL^{-2}\), \([V]=L^{-1}\), \([O]=N_oL^{-2}\), and
\([\Pi]=P_a\).  The prescribed tumor profile \(\Theta\) is dimensionless.
The hierarchy
\begin{equation}
C\quad\longrightarrow\quad V\quad\longrightarrow\quad
\Phi[V]=V\widehat G(Q[V]),\qquad 0\leq\Phi[V]\leq V,
\label{eq:levels}
\end{equation}
distinguishes motile endothelial occupancy, formed vascular structure, and pressure-supported functional structure.
Here, $\Phi[V]$ represents the flow-functional fraction of the vascular structure.
The variable \(P\) is a phenomenological net degradation
closure rather than an asserted fast-binding reduction.  Figure~\ref{fig:architecture}
shows how the three levels couple through the nonlocal perfusion closure
analyzed in Section~\ref{sec:perfusion}.

For a radial \(\widehat J\in C_c^2(\R^2)\), nonnegative, positive near the
origin, and of unit integral, set
\begin{equation*}
\widehat{\mathcal J}_{\ell_d}V(X)=
\frac{\int_{\widehat\Omega}\ell_d^{-2}
\widehat J((X-Y)/\ell_d)V(Y)\dd Y}
{\int_{\widehat\Omega}\ell_d^{-2}
\widehat J((X-Y)/\ell_d)\dd Y}.
\end{equation*}
With conductivity unit \(L^2(P_aT)^{-1}\), define
\begin{align*}
\widehat K[V]&=\widehat K_b+\widehat K_v
\left(\widehat{\mathcal J}_{\ell_d}V/V_0\right)^m,
&-\nabla_X\!\cdot(\widehat K[V]\nabla_X\Pi)&=0,\\
Q[V]&=\widehat K_v
\left(\widehat{\mathcal J}_{\ell_d}V/V_0\right)^m
\abs{\nabla_X\Pi[V]},
&\widehat G(Q)&=\frac{Q}{Q+Q_c}.
\end{align*}
Pressure equals \(\Pi_A>\Pi_V\) on effective arterial and venous boundary
parts and has zero conormal flux elsewhere.  The background conductivity is
excluded from \(Q\), so it regularizes pressure without being counted as
blood flow.
Let
\begin{equation*}
\begin{aligned}
\widehat M(C)=C(1-C/C_0),\quad
&\widehat S(A)=\widehat\kappa_s^{-1}\log(1+\widehat\kappa_sA),\quad
\widehat H_P(O)=\frac{O}{\widehat K_P+O},\\
&\widehat H_H(O)=\frac{(\widehat O_h-O)_+}
{\widehat K_H+(\widehat O_h-O)_+}.
\end{aligned}
\end{equation*}
The dimensional biological balances are
\begin{align*}
C_{\mathsf t}={}&\nabla_X\!\cdot[\widehat D_C\nabla_XC
-\widehat X_A\widehat M(C)\nabla_X\widehat S(A)
-\widehat X_F\widehat M(C)\nabla_XF] \nonumber \\
{}&+\widehat r_C\widehat M(C)\widehat H_P(O)
-\widehat\eta C(1-V/V_0),\\
P_{\mathsf t}={}&\widehat D_P\Delta_XP
+\widehat\Sigma_P\Theta(\theta_0+\theta_1C/C_0)-\widehat\mu_PP,\\
F_{\mathsf t}={}&\widehat\rho_F(F_0-F)-\widehat k_FPF,\\
\widehat\tau_AA_{\mathsf t}={}&\widehat L_A^2\Delta_XA
+\widehat\Sigma_A\Theta\widehat H_H(O)-d_AA-k_ACA,\\
V_{\mathsf t}={}&Y_V\widehat\eta C(1-V/V_0)
-\widehat\delta_V[1-\widehat G(Q[V])]V,\\
\widehat\tau_OO_{\mathsf t}={}&\widehat L_O^2\Delta_XO
-\widehat\Gamma_T\Theta\frac{O}{\widehat K_T+O}
-\widehat\Gamma_C\frac{C}{C_0}\frac{O}{\widehat K_C+O}
+\widehat\Lambda\Phi[V](\widehat O_b-O).
\end{align*}
All terms in each balance have the same units as its left-hand side.  In
particular, \([\widehat X_A]=L^2[A]^{-1}T^{-1}\),
\([\widehat X_F]=L^2[F]^{-1}T^{-1}\),
\([\widehat k_F]=[P]^{-1}T^{-1}\), and
\([\widehat\Lambda]=[V]^{-1}\).  The endothelial flux and the diffusive
fluxes of \(P,A,O\) have zero normal component. Table~\ref{tab:variables-units} lists the dimensional variables and their reference scales, and Table~\ref{tab:dimensional-parameters} lists the dimensional motility, protease, matrix, vessel, signal, oxygen, and pressure parameters.

\begin{table}[htbp]
\caption{Dimensional variables and reference scales.}\label{tab:variables-units}
\begin{tabular}{@{}p{0.18\textwidth}p{0.51\textwidth}p{0.23\textwidth}@{}}
\toprule
Symbol & Meaning & Unit\\
\midrule
\(X,L_0,\ell_d\) & position, tissue scale, averaging length & \(L\)\\
\(\mathsf t,t_0\) & dimensional and reference time & \(T\)\\
\(C,C_0\) & endothelial density and saturation scale & \(N_cL^{-2}\)\\
\(P,P_0\) & net proteolytic activity and scale & \(N_pL^{-2}\)\\
\(F,F_0\) & extracellular matrix (ECM) density and intact-matrix scale & \(M_fL^{-2}\)\\
\(A,A_0\) & VEGF-like concentration and scale & \(N_aL^{-2}\)\\
\(V,V_0\) & vessel length per tissue area and scale & \(L^{-1}\)\\
\(O,O_0\) & oxygen concentration and scale & \(N_oL^{-2}\)\\
\(\Pi,\Pi_0\) & vascular pressure and scale & \(P_a\)\\
\(Q,Q_0\) & Darcy-scale flow magnitude and scale & \(LT^{-1}\)\\
\(\Theta\) & stationary tumor profile & \(1\)\\
\botrule
\end{tabular}
\end{table}

\begin{table}[htbp]
\caption{Dimensional motility, protease, matrix, vessel, signal, oxygen, and pressure parameters.}\label{tab:dimensional-parameters}
\begin{tabular}{@{}p{0.21\textwidth}p{0.49\textwidth}p{0.22\textwidth}@{}}
\toprule
Symbol & Meaning & Unit\\
\midrule
\(\widehat D_C,\widehat D_P\) & endothelial and protease diffusivities & \(L^2T^{-1}\)\\
\(\widehat X_A\) & sensitivity to \(\widehat S(A)\) & \(L^2[A]^{-1}T^{-1}\)\\
\(\widehat X_F\) & sensitivity to ECM & \(L^2[F]^{-1}T^{-1}\)\\
\(\widehat\kappa_s\) & inverse VEGF sensing scale & \([A]^{-1}\)\\
\(\widehat r_C,\widehat\eta\) & proliferation and lumenization rates & \(T^{-1}\)\\
\(\widehat K_P\) & proliferation oxygen half-saturation & \([O]\)\\
\(\widehat\Sigma_P\) & protease source scale & \([P]T^{-1}\)\\
\(\theta_0,\theta_1\) & source weights & \(1\)\\
\(\widehat\mu_P,\widehat\rho_F,\widehat\delta_V\) & decay, renewal, and regression rates & \(T^{-1}\)\\
\(\widehat k_F\) & proteolytic degradation coefficient & \([P]^{-1}T^{-1}\)\\
\(Y_V\) & lumenized-density yield & \([V][C]^{-1}\)\\
\(\widehat\tau_A,\widehat\tau_O\) & VEGF and oxygen relaxation times & \(T\)\\
\(\widehat L_A^2,\widehat L_O^2\) & relaxation-scaled diffusion & \(L^2\)\\
\(\widehat\Sigma_A\) & VEGF production per relaxation & \([A]\)\\
\(d_A,k_A\) & VEGF decay and endothelial uptake & \(1,[C]^{-1}\)\\
\(\widehat O_h,\widehat O_b,\widehat K_H\) & hypoxia, blood, and response scales & \([O]\)\\
\(\widehat\Gamma_T,\widehat\Gamma_C\) & oxygen consumption per relaxation & \([O]\)\\
\(\widehat K_T,\widehat K_C\) & consumption half-saturations & \([O]\)\\
\(\widehat\Lambda\) & vessel-to-tissue exchange coefficient & \([V]^{-1}\)\\
\(\widehat K_b,\widehat K_v\) & background and vessel conductivities & \(L^2(P_aT)^{-1}\)\\
\(Q_c\) & flow-gating threshold & \(LT^{-1}\)\\
\(\Pi_A,\Pi_V\) & imposed boundary pressures & \(P_a\)\\
\(m\) & conductivity-density exponent & \(1\)\\
\(\widehat J,\widehat J_{\ell_d}\) & reference and scaled kernels & \(1,L^{-2}\)\\
\botrule
\end{tabular}
\end{table}

\subsection{Scaling and nondimensional system}
Choose scales \(L_0,t_0,C_0,P_0,F_0,A_0,V_0,O_0,\Pi_0\), set
\begin{equation*}
x=X/L_0,\quad t=\mathsf t/t_0,\quad
(c,p,f,a,v,o,\pi)=
(C/C_0,P/P_0,F/F_0,A/A_0,V/V_0,O/O_0,\Pi/\Pi_0),
\end{equation*}
and choose \(V_0=Y_VC_0\).  With a conductivity scale
\(\widehat K_{\rm ref}\), put
\(Q_0=\widehat K_{\rm ref}\Pi_0/L_0\) and \(\ell=\ell_d/L_0\).
The following dimensionless coefficients provide the parameterization used below:
\begin{align*}
D_c&=t_0\widehat D_C/L_0^2,&
\chi_a&=t_0\widehat X_AA_0/L_0^2,&
\chi_f&=t_0\widehat X_FF_0/L_0^2,&
\kappa_s&=\widehat\kappa_sA_0,\\
r_c&=t_0\widehat r_C,&K_P&=\widehat K_P/O_0,&
\eta&=t_0\widehat\eta,&D_p&=t_0\widehat D_P/L_0^2,\\
\sigma_p&=t_0\widehat\Sigma_P/P_0,&
\mu_p&=t_0\widehat\mu_P,&\rho_f&=t_0\widehat\rho_F,&
\kappa_f&=t_0\widehat k_FP_0,\\
\tau_a&=\widehat\tau_A/t_0,&D_a&=\widehat L_A^2/L_0^2,&
\sigma_a&=\widehat\Sigma_A/A_0,&
(\mu_a,\kappa_a)&=(d_A,k_AC_0),\\
\delta_v&=t_0\widehat\delta_V,&
(o_h,K_H)&=(\widehat O_h,\widehat K_H)/O_0,&
\tau_o&=\widehat\tau_O/t_0,&D_o&=\widehat L_O^2/L_0^2,\\
(\gamma_T,\gamma_c)&=(\widehat\Gamma_T,\widehat\Gamma_C)/O_0,&
(K_T,K_c)&=(\widehat K_T,\widehat K_C)/O_0,&
\lambda&=\widehat\Lambda V_0,&o_b&=\widehat O_b/O_0,\\
\eps&=\widehat K_b/\widehat K_{\rm ref},&
K_v&=\widehat K_v/\widehat K_{\rm ref},&
q_*&=Q_c/Q_0.&
\end{align*}
The coefficients above provide a convenient dimensionless
parameterization under the stated choice of reference scales; they are
not intended to constitute a minimal set of dimensionless invariants.
For the pressure--perfusion subsystem, two elementary invariances reduce
the number of effective hydrodynamic parameters. First, the pressure is
defined only up to an additive constant, so the pressure field enters
the flow law through the pressure drop
\(\Delta\pi:=\pi_A-\pi_V>0\). Second, multiplication of the conductivity
coefficient \(K[v]\) by a common positive factor leaves the pressure
equation unchanged. Consequently, after factoring out the common
conductivity scale, the pressure--perfusion closure depends on
\[
\frac{\varepsilon}{K_v},
\qquad
\frac{q_*}{K_v(\pi_A-\pi_V)},
\qquad
\ell,
\qquad
m,
\]
together with the scaled domain, the pressure-boundary partition, and
the shape of the kernel \(J\).

In the protease equation, the source coefficients enter only through
the products \(\sigma_p\theta_0\) and \(\sigma_p\theta_1\). Moreover,
each of the VEGF and oxygen equations is homogeneous with respect to a
common positive multiplication of its five displayed coefficients.
These are algebraic invariances of the chosen dimensionless
parameterization. Accordingly, the displayed parameter
list is explicit but not minimal.
For
\begin{equation*}
\begin{aligned}
&M(c)=c(1-c),\quad S(a)=\kappa_s^{-1}\log(1+\kappa_sa),\quad
H_P(o)=\frac{o}{K_P+o},\\
&H_H(o)=\frac{(o_h-o)_+}{K_H+(o_h-o)_+},\quad
G(q)=\frac{q}{q+q_*},
\end{aligned}
\end{equation*}
define the boundary-normalized average
\begin{equation}
\Jell v(x)=\frac{\int_\Omega J_\ell(x-y)v(y)\dd y}
{\int_\Omega J_\ell(x-y)\dd y},\qquad
J_\ell(z)=\ell^{-2}J(z/\ell),
\label{eq:average}
\end{equation}
and
\begin{align*}
K[v]&=\eps+K_v(\Jell v)^m,&
-\nabla\!\cdot(K[v]\nabla\pi[v])&=0,\\
q[v]&=K_v(\Jell v)^m\abs{\nabla\pi[v]},&
\Ffun[v]&=vG(q[v]).
\end{align*}
The pressure traces are \(\pi_A>\pi_V\) on \(\Gamma_A,\Gamma_V\), with
zero conormal flux on \(\Gamma_N\).  The biological system is
\begin{align}
c_t={}&\nabla\!\cdot[D_c\nabla c-\chi_aM(c)\nabla S(a)
-\chi_fM(c)\nabla f]+r_cM(c)H_P(o)-\eta c(1-v),
\label{eq:c}\\
p_t={}&D_p\Delta p+\sigma_pT(\theta_0+\theta_1c)-\mu_pp,
\label{eq:p}\\
f_t={}&\rho_f(1-f)-\kappa_fpf,
\label{eq:f}\\
\tau_aa_t={}&D_a\Delta a+\sigma_aTH_H(o)-\mu_aa-\kappa_aca,
\label{eq:a}\\
v_t={}&\eta c(1-v)-\delta_v[1-G(q[v])]v,
\label{eq:v}\\
\tau_oo_t={}&D_o\Delta o-\gamma_TT\frac{o}{K_T+o}
-\gamma_cc\frac{o}{K_c+o}+\lambda\Ffun[v](o_b-o).
\label{eq:o}
\end{align}
The complete conservative flux in \eqref{eq:c}, and the diffusive fluxes in
\eqref{eq:p}, \eqref{eq:a}, and \eqref{eq:o}, have zero normal component.
\subsection{Standing assumptions and scope}
The following assumptions are used throughout unless a theorem states a
smaller subset.

\begin{enumerate}
\item[\textbf{(A1)}] \(\Omega\subset\R^2\) is bounded, connected, and
\(C^{1,1}\).  The closures of pairwise disjoint relatively open sets
\(\Gamma_A,\Gamma_V,\Gamma_N\) cover its boundary, and
\(\Gamma_A,\Gamma_V\) have positive one-dimensional Hausdorff (surface)
measure.
\item[\textbf{(A2)}] The pressure data admit an \(H^1(\Omega)\) lifting.
All pressure and biological boundary conditions are those stated above.
\item[\textbf{(A3)}] \(J\in C_c^2(\R^2)\) is radial, nonnegative,
positive near zero, and has unit integral; \(\eps,\ell,K_v,q_*>0\) are fixed
and \(m\geq1\).
\item[\textbf{(A4)}] All coefficients are nonnegative.  The four diffusion
coefficients, $\tau_a$, $\tau_o$, $\mu_p$, $\rho_f$, $\mu_a$, $\delta_v$, and every
constant occurring in a denominator are positive.
\item[\textbf{(A5)}] \(T\in L^{\infty}(\Omega)\), \(0\leq T\leq1\),
\(|\{T>0\}|>0\), and \(\gamma_T>0\).
\item[\textbf{(A6)}] The initial data satisfy
\(f_0\in H^1(\Omega)\),
\(0\leq c_0,f_0,v_0\leq1\), \(0\leq o_0\leq o_b\), and
\(p_0,a_0\in L^\infty(\Omega)\) with \(p_0,a_0\geq0\).
\item[\textbf{(A7)}] Approximation uses the globally Lipschitz extension
\(M_*(z)=M(z)\) on \([0,1]\) and \(M_*(z)=0\) outside this interval.  It
disappears after the invariant-region estimate.
The model treats the tumor as prescribed and stationary.  It excludes
discrete anastomosis, exact connectivity, intravascular oxygen advection,
red-cell rheology, individual radii, tumor growth, mechanics, treatment, and
the limit \(\eps\to0\).  Logarithmic sensing reduces marginal sensitivity at
large \(a\), but does not bound taxis independently of \(\nabla a\).
\end{enumerate}

\section{The nonlocal pressure--perfusion map}
\label{sec:perfusion}
Write
\begin{equation*}
d_\ell(x)=\int_\Omega J_\ell(x-y)\dd y,
\qquad N_\ell h(x)=\int_\Omega J_\ell(x-y)h(y)\dd y,
\qquad \Jell h=N_\ell h/d_\ell.
\end{equation*}

A bounded \(C^{1,1}\) domain has constants \(r_\Omega,c_\Omega>0\) such
that \(|\Omega\cap B_r(x)|\geq c_\Omega r^2\) for
\(x\in\overline\Omega\) and \(r\leq r_\Omega\).  Choose
\(\rho_J,j_J>0\) with \(J\geq j_J\) on \(B_{\rho_J}\), and set
\begin{align*}
r_\ell&=\min\{r_\Omega,\rho_J\ell\},&
\underline d_\ell&=j_Jc_\Omega\ell^{-2}r_\ell^2,\\
A_{\ell,0}&=\frac{\norm J_2}{\underline d_\ell\ell},&
A_{\ell,1}&=\frac{\norm{\nabla J}_2}{\underline d_\ell\ell^2}
+\frac{\norm J_2\norm{\nabla J}_1}
{\underline d_\ell^2\ell^2}.
\end{align*}
\begin{lemma}[Normalized kernel]
\label{lem:kernel}
Under (A1) and (A3),
\begin{equation*}
0<\underline d_\ell\leq d_\ell\leq1,
\qquad
\norm{\Jell h}_{\infty}\leq A_{\ell,0}\norm h_2,
\qquad
\norm{\nabla\Jell h}_{\infty}\leq A_{\ell,1}\norm h_2.
\end{equation*}
Consequently, \(\Jell:L^2(\Omega)\to W^{1,\infty}(\Omega)\) is
Lipschitz, and \(0\leq\Jell v\leq1\) whenever \(0\leq v\leq1\).
\end{lemma}

\begin{proof}
By (A3), \(J\in C_c^2(\mathbb R^2)\), \(J\ge0\), and \(J\) is positive in a neighborhood of the origin. Hence there exist \(\rho_J,j_J>0\) such that $J(z)\ge j_J$ for $|z|\le\rho_J$.
By the uniform interior-density property of the bounded \(C^{1,1}\) domain \(\Omega\), there exist \(r_\Omega,c_\Omega>0\) such that
\[
|\Omega\cap B_r(x)|\ge c_\Omega r^2,
\qquad x\in\overline\Omega,\quad 0<r\le r_\Omega.
\]
Set $r_\ell=\min\{r_\Omega,\rho_J\ell\}$.
If \(y\in\Omega\cap B_{r_\ell}(x)\), then
$\left|\dfrac{x-y}{\ell}\right|
\le\dfrac{r_\ell}{\ell}\le\rho_J$,
and hence
$J_\ell(x-y)=\ell^{-2}J((x-y)/\ell)\ge j_J\ell^{-2}$.
Therefore
$$
d_\ell(x)
=\int_\Omega J_\ell(x-y)\,\dd y
\ge
j_J\ell^{-2}|\Omega\cap B_{r_\ell}(x)|
\ge
j_Jc_\Omega\ell^{-2}r_\ell^2
=\underline d_\ell.
$$
Moreover, since \(J_\ell\ge0\) and \(\int_{\mathbb R^2}J_\ell=1\),
$d_\ell(x)\le\int_{\mathbb R^2}J_\ell(z)\,\dd z=1$, establishing the first estimate.

Next, the scaling \(J_\ell(z)=\ell^{-2}J(z/\ell)\) gives
\[
\|J_\ell\|_2=\ell^{-1}\|J\|_2,\qquad
\|\nabla J_\ell\|_2=\ell^{-2}\|\nabla J\|_2,\qquad
\|\nabla J_\ell\|_1=\ell^{-1}\|\nabla J\|_1.
\]
For \(h\in L^2(\Omega)\), Cauchy--Schwarz yields, for every \(x\in\overline\Omega\), $|N_\ell h(x)|
\le
\|J_\ell\|_2\|h\|_2$,
and, differentiating under the integral sign,
$$
\nabla N_\ell h(x)
=
\int_\Omega \nabla J_\ell(x-y)h(y)\,\dd y,
$$
so $|\nabla N_\ell h(x)|
\le
\|\nabla J_\ell\|_2\|h\|_2$.
Hence $\|N_\ell h\|_\infty\le\|J_\ell\|_2\|h\|_2$ and $\|\nabla N_\ell h\|_\infty
\le\|\nabla J_\ell\|_2\|h\|_2$. 
Furthermore,
$$
\nabla d_\ell(x)
=
\int_\Omega \nabla J_\ell(x-y)\,\dd y,
$$
and therefore $\|\nabla d_\ell\|_\infty
\le
\|\nabla J_\ell\|_1$.

Since \(d_\ell\ge\underline d_\ell>0\), $\Jell h=\dfrac{N_\ell h}{d_\ell}$
and $\nabla\Jell h=\dfrac{\nabla N_\ell h}{d_\ell}-\dfrac{N_\ell h\,\nabla d_\ell}{d_\ell^2}$.
Consequently,
\[
\|\Jell h\|_\infty\le
\frac{\|J_\ell\|_2}{\underline d_\ell}\|h\|_2=
\frac{\|J\|_2}{\underline d_\ell\ell}\|h\|_2=
A_{\ell,0}\|h\|_2,
\]
and
\[
\|\nabla\Jell h\|_\infty
\le
\left(
\frac{\|\nabla J_\ell\|_2}{\underline d_\ell}
+
\frac{\|J_\ell\|_2\|\nabla J_\ell\|_1}
{\underline d_\ell^2}
\right)\|h\|_2=
\left(
\frac{\|\nabla J\|_2}{\underline d_\ell\ell^2}
+
\frac{\|J\|_2\|\nabla J\|_1}
{\underline d_\ell^2\ell^2}
\right)\|h\|_2=
A_{\ell,1}\|h\|_2.
\]
Thus \(\Jell:L^2(\Omega)\to W^{1,\infty}(\Omega)\) is a bounded linear operator, hence Lipschitz.

Finally, if \(0\le v\le1\) a.e., then \(J_\ell\ge0\) implies $0\le N_\ell v\le N_\ell1=d_\ell$.
Since \(d_\ell>0\),
$0\le\Jell v=\dfrac{N_\ell v}{d_\ell}\le1$.
This proves the lemma.
\end{proof}

Let
\begin{equation*}
\mathcal A=\{v\in L^2(\Omega):0\leq v\leq1\ \text{a.e.}\},
\qquad
\mathcal V=\{\varphi\in H^1(\Omega):\operatorname{Tr}\varphi=0
\text{ on }\Gamma_A\cup\Gamma_V\}.
\end{equation*}
Fix an \(H^1\) lifting \(\pi_D\).  If \(C_P\) is the Poincar\'e
constant on \(\mathcal V\), put \(C_E=(1+C_P^2)^{1/2}\),
\(\overline K=\eps+K_v\), and
\begin{align*}
M_\nabla&=\left(1+\frac{\overline K}{\eps}\right)
\norm{\nabla\pi_D}_2,&
M_H&=\norm{\pi_D}_{H^1}+C_E\frac{\overline K}{\eps}
\norm{\nabla\pi_D}_2,\\
L_K&=K_vmA_{\ell,0},&
L_\nabla&=\frac{L_KM_\nabla}{\eps},\qquad
L_\pi=C_EL_\nabla,\\
L_q&=K_v(mA_{\ell,0}M_\nabla+L_\nabla),&
L_{\Ffun}&=1+L_q/q_*.
\end{align*}
The next result is the interface used by both compactness arguments.

\begin{theorem}[Pressure-to-perfusion map]
\label{thm:perfusion}
Under (A1)--(A3), every \(v\in\mathcal A\) determines a unique
\(\pi[v]\in\pi_D+\mathcal V\) satisfying
\begin{equation}
\int_\Omega K[v]\nabla\pi[v]\cdot\nabla\varphi\dd x=0
\qquad(\varphi\in\mathcal V).
\label{eq:weak-pressure}
\end{equation}
It obeys
\begin{equation}
\norm{\pi[v]}_{H^1}\leq M_H,\qquad
\norm{q[v]}_2\leq K_vM_\nabla.
\label{eq:map-bounds}
\end{equation}
For \(v_1,v_2\in\mathcal A\),
\begin{multline*}
\norm{\pi[v_1]-\pi[v_2]}_{H^1}
+\norm{q[v_1]-q[v_2]}_2
+\norm{\Ffun[v_1]-\Ffun[v_2]}_2\\
\leq (L_\pi+L_q+L_{\Ffun})\norm{v_1-v_2}_2.
\end{multline*}
Moreover, \(0\leq\Ffun[v]\leq v\leq1\).  If
\(t\mapsto v(t)\) is strongly measurable in \(L^2\), then the associated
\(\pi,q,\Ffun\) are strongly measurable in \(H^1,L^2,L^2\), respectively.
\end{theorem}

\begin{proof}
Lemma~\ref{lem:kernel} yields
\(\eps\leq K[v]\leq\overline K\).  Writing \(\pi=\pi_D+w\), the equation
for \(w\in\mathcal V\) is
\[
\int_\Omega K[v]\nabla w\cdot\nabla\varphi
=-\int_\Omega K[v]\nabla\pi_D\cdot\nabla\varphi.
\]
Lax--Milgram applies in the gradient norm.  Testing by \(w\) gives
\(\norm{\nabla w}_2\leq\overline K\norm{\nabla\pi_D}_2/\eps\),
which proves \eqref{eq:map-bounds}; the natural boundary condition is the
zero conormal flux on \(\Gamma_N\).
For \(a_i=(\Jell v_i)^m\), \(K_i=\eps+K_va_i\), the Lipschitz property of
\(s^m\) on \([0,1]\) gives
\begin{equation*}
\norm{a_1-a_2}_\infty\leq mA_{\ell,0}\norm{v_1-v_2}_2,
\qquad
\norm{K_1-K_2}_\infty\leq L_K\norm{v_1-v_2}_2.
\end{equation*}
Subtracting the pressure equations and testing by
\(u=\pi[v_1]-\pi[v_2]\in\mathcal V\) yields
\[
\norm{\nabla u}_2\leq L_\nabla\norm{v_1-v_2}_2,
\qquad
\norm u_{H^1}\leq L_\pi\norm{v_1-v_2}_2.
\]
Since \(\lvert\abs\xi-\abs\zeta\rvert\leq\abs{\xi-\zeta}\),
\[
\norm{q[v_1]-q[v_2]}_2
\leq K_v\bigl(\norm{a_1-a_2}_\infty M_\nabla
+\norm{\nabla u}_2\bigr)
\leq L_q\norm{v_1-v_2}_2.
\]
Finally, \(0\leq G\leq1\), \(\operatorname{Lip}(G)=q_*^{-1}\), and
\begin{equation*}
\Ffun[v_1]-\Ffun[v_2]=(v_1-v_2)G(q[v_1])
+v_2\{G(q[v_1])-G(q[v_2])\}
\end{equation*}
give the last estimate and the order interval.  The admissible set is closed
in \(L^2\); the three maps are Lipschitz into separable spaces, so composition
with a strongly measurable path proves measurability.
\end{proof}

The constants deteriorate when \(\eps\downarrow0\) or \(\ell\downarrow0\):
for small \(\ell\), \(A_{\ell,0}=O(\ell^{-1})\) and
\(A_{\ell,1}=O(\ell^{-2})\), while \(M_\nabla=O(\eps^{-1})\), and the
displayed stability constant contains a term of order \(\eps^{-2}\), namely,
\(L_\nabla=L_KM_\nabla/\eps\).  Neither
singular limit is considered.  The theorem proves continuity, not
compactness, of \(v\mapsto\Ffun[v]\), because the factor \(v\) is not
smoothed.

\section{Global bounded weak solutions}
\label{sec:global}
Fix \(T_f>0\), write \(Q_T=\Omega\times(0,T_f)\), and assume
(A1)--(A7).  A weak solution has
\begin{align}
c,p,a,o&\in L^2(0,T_f;H^1)\cap L^\infty(Q_T),&
c_t,p_t,a_t,o_t&\in L^2(0,T_f;(H^1)^*),
\label{eq:weak-reg-diffusive}\\
f&\in L^\infty(0,T_f;H^1)\cap W^{1,\infty}(0,T_f;L^2),&
v&\in W^{1,\infty}(0,T_f;L^2)\cap L^\infty(Q_T),
\label{eq:weak-reg-ode}\\
\pi&\in L^\infty(0,T_f;H^1).&&
\label{eq:weak-reg-pressure}
\end{align}
For a.e. \(t\), every \(\varphi,\psi,\xi,\zeta\in H^1(\Omega)\),
and \(h\in L^2(\Omega)\), it satisfies
\begin{align}
\langle c_t,\varphi\rangle&+\int_\Omega[D_c\nabla c-\chi_aM(c)\nabla S(a)
-\chi_fM(c)\nabla f]\cdot\nabla\varphi=\int_\Omega[r_cM(c)H_P(o)-\eta c(1-v)]\varphi,
\label{eq:weak-c-compact}\\
\langle p_t,\psi\rangle&+D_p(\nabla p,\nabla\psi)
=(\sigma_pT(\theta_0+\theta_1c)-\mu_pp,\psi),
\label{eq:weak-p-compact}\\
(f_t,h)&=(\rho_f(1-f)-\kappa_fpf,h),
\label{eq:weak-f-compact}\\
\tau_a\langle a_t,\xi\rangle&+D_a(\nabla a,\nabla\xi)
+((\mu_a+\kappa_ac)a,\xi)
=(\sigma_aTH_H(o),\xi),
\label{eq:weak-a-compact}\\
(v_t,h)&=(\eta c(1-v)-\delta_v[1-G(q[v])]v,h),
\label{eq:weak-v-compact}\\
\tau_o\langle o_t,\zeta\rangle&+D_o(\nabla o,\nabla\zeta)
+\left(\gamma_TT\frac{o}{K_T+o}
+\gamma_cc\frac{o}{K_c+o},\zeta\right)=(\lambda\Ffun[v](o_b-o),\zeta).
\label{eq:weak-o-compact}
\end{align}
Here and below unlabeled integrals are over \(\Omega\), \(\pi=\pi[v]\)
solves \eqref{eq:weak-pressure}, and the initial values are attained in
\(L^2\).
The following result shows that the nondiffusive vascular component does not
destroy global compactness.
\begin{theorem}[Global bounded weak solutions]
\label{thm:global}
For fixed \(\tau_a,\tau_o,\eps,\ell>0\), assumptions (A1)--(A7) imply
existence of a weak solution on every finite interval with
\begin{equation}
0\leq c,f,v\leq1,\qquad 0\leq o\leq o_b,\qquad p,a\geq0,
\label{eq:invariant}
\end{equation}
and
\begin{equation}
\norm p_\infty\leq
\max\left\{\norm{p_0}_\infty,
\frac{\sigma_p(\theta_0+\theta_1)}{\mu_p}\right\},\qquad
\norm a_\infty\leq
\max\left\{\norm{a_0}_\infty,\frac{\sigma_a}{\mu_a}\right\}.
\label{eq:comparison-bounds}
\end{equation}
The solution has the regularity stated in \eqref{eq:weak-reg-diffusive}--\eqref{eq:weak-reg-pressure}; in addition,
\(q[v]\in L^\infty(0,T_f;L^2)\).  Uniqueness of weak solutions is not
asserted.
\end{theorem}
The full proof is provided in Appendix~\ref{app:approximation}. We then record a compactness mechanism for \(v\) that will used in the proof of Theorem~\ref{thm:global}.

\begin{lemma}[Vascular stability]
\label{lem:vascular}
Let \(0\leq b_i,v_i\leq1\), with
\[
(v_i)_t=\eta b_i(1-v_i)-\delta_v[1-G(q[v_i])]v_i.
\]
Then, for a.e. time,
\begin{equation}
\frac{\dd}{\dd t}\norm{v_1-v_2}_2^2
\leq \eta\norm{b_1-b_2}_2^2
+\left(\eta+\frac{2\delta_vL_q}{q_*}\right)
\norm{v_1-v_2}_2^2.
\label{eq:vascular-stability-compact}
\end{equation}
Consequently, strong convergence of \(b_n\) in \(L^2(Q_T)\), together
with strong convergence of the initial vessels, implies
\(v_n\to v\) in \(C([0,T_f];L^2)\).
\end{lemma}
\begin{proof}
Put \(z=v_1-v_2\), \(q_i=q[v_i]\).  The exact decompositions
\begin{align*}
b_1(1-v_1)-b_2(1-v_2)&=(b_1-b_2)(1-v_1)-b_2z,\\
[1-G(q_1)]v_1-[1-G(q_2)]v_2
&=[1-G(q_1)]z+v_2[G(q_2)-G(q_1)]
\end{align*}
are tested by \(z\).  Discard the two nonpositive quadratic terms, use
\(0\leq b_i,v_i\leq1\),
\(\operatorname{Lip}(G)=q_*^{-1}\), Theorem~\ref{thm:perfusion}, and
Young's inequality.  This gives \eqref{eq:vascular-stability-compact};
Gr\"onwall gives the last assertion.
\end{proof}

In the applications below, the comparison solution is constructed through
a globally Lipschitz extension of the vessel equation.  Let
\(\mathsf C(s)=\min\{1,\max\{0,s\}\}\), acting pointwise, and, for a
strongly measurable driver \(0\leq b\leq1\), define
\begin{equation}
\mathcal B_b(w)=\eta b[1-\mathsf C(w)]
-\delta_v\{\mathsf C(w)-\Ffun[\mathsf C(w)]\},
\qquad w\in L^2(\Omega).
\label{eq:vessel-extension}
\end{equation}
This removes the domain mismatch between the vessel equation and
Theorem~\ref{thm:perfusion}: \(\mathcal B_b\) is defined on all of
\(L^2(\Omega)\), while \(\Ffun\) is evaluated only at
\(\mathsf C(w)\in\mathcal A\), and \(\mathcal B_b(w)\) reproduces the
right-hand side of \eqref{eq:v} whenever \(w\in\mathcal A\).  Since
\(\mathsf C:L^2(\Omega)\to\mathcal A\) is \(1\)-Lipschitz and
\(\Ffun:\mathcal A\to L^2(\Omega)\) is \(L_{\Ffun}\)-Lipschitz by
Theorem~\ref{thm:perfusion},
\begin{equation}
\norm{\mathcal B_b(w_1)-\mathcal B_b(w_2)}_2
\leq\{\eta+\delta_v(1+L_{\Ffun})\}\norm{w_1-w_2}_2.
\label{eq:vessel-extension-lipschitz}
\end{equation}
For each fixed \(w\in L^2(\Omega)\), the map
\(t\mapsto\mathcal B_{b(t)}(w)\) is strongly measurable.  By
\eqref{eq:vessel-extension-lipschitz},
\[
\norm{\mathcal B_{b(t)}(w_1)-\mathcal B_{b(t)}(w_2)}_2
\le L_B\norm{w_1-w_2}_2,
\qquad
L_B:=\eta+\delta_v(1+L_{\Ffun}),
\]
uniformly in \(t\).  Furthermore, since
\(0\le b(t),\mathsf C(w),\Ffun[\mathsf C(w)]\le1\),
\[
\norm{\mathcal B_{b(t)}(w)}_2
\le (\eta+\delta_v)|\Omega|^{1/2}
\]
for all \(t\) and \(w\).  Thus
\(\mathcal B\) satisfies the Carath\'eodory hypotheses on the Banach
space \(L^2(\Omega)\), with a global Lipschitz constant and a global
growth bound.  Hence the nonautonomous Cauchy problem
\[
v_t=\mathcal B_{b(t)}(v),\qquad v(0)=v_0\in\mathcal A,
\]
admits a unique global solution.  Since
\(v_t=\mathcal B_{b(t)}(v(t))\) is uniformly bounded in \(L^2(\Omega)\),
the solution satisfies $v\in W^{1,\infty}(0,T_f;L^2(\Omega))$.
Moreover, the interval \([0,1]\) is invariant for the extended equation.
Indeed, on \(\{v<0\}\) one has \(\mathsf C(v)=0\) and
\(\Ffun[\mathsf C(v)]=\Ffun[0]=0\), since
\(\Ffun[w]=wG(q[w])\). Hence $\mathcal B_b(v)=\eta b\geq0$ on $\{v<0\}$.
On \(\{v>1\}\), one has \(\mathsf C(v)=1\), and
\(0\leq\Ffun[1]\leq1\), so that $\mathcal B_b(v)
=-\delta_v\{1-\Ffun[1]\}\leq0$ on $\{v>1\}$.
Let
\[
v^-:=\max\{-v,0\},
\qquad
(v-1)^+:=\max\{v-1,0\}.
\]
Testing the extended equation with \(-v^-\) and \((v-1)^+\),
respectively, and using the standard chain rule for Lipschitz
truncations, gives
\[
\frac12\frac{\dd}{\dd t}\|v^-\|_2^2\leq0,
\qquad
\frac12\frac{\dd}{\dd t}\|(v-1)^+\|_2^2\leq0.
\]
Since \(0\leq v_0\leq1\), we have
\(v_0^-=(v_0-1)^+=0\), and consequently $v^-=0$, $(v-1)^+=0$, hence
$0\leq v\leq1$ a.e. in $Q_T$.
Therefore \(\mathsf C(v)=v\) a.e., and the extended equation coincides
with the original vessel ODE \eqref{eq:v}. Conversely, every
\(\mathcal A\)-valued solution of \eqref{eq:v} solves the extended
equation, so uniqueness transfers. In particular, the limit field \(v\)
in Theorems~\ref{thm:global} and~\ref{thm:fast} is the unique
\(\mathcal A\)-valued solution of the vessel ODE driven by \(c\).

\section{Simultaneous fast oxygen--VEGF limit}
\label{sec:fast}
The limiting signal equations, with homogeneous Neumann conditions, are
\begin{align}
-D_o\Delta o+\gamma_TT\frac{o}{K_T+o}
+\gamma_cc\frac{o}{K_c+o}+\lambda\Ffun[v]o
&=\lambda\Ffun[v]o_b,
\label{eq:elliptic-o-compact}\\
-D_a\Delta a+(\mu_a+\kappa_ac)a&=\sigma_aTH_H(o).
\label{eq:elliptic-a-compact}
\end{align}
Their solvability does not assume positive perfusion everywhere.

\begin{lemma}[Weighted Poincar\'e inequality]
\label{lem:weighted-poincare}
Under (A1) and (A5),
\begin{equation}
\norm z_2^2\leq C_{\mathrm{wp}}\left(\norm{\nabla z}_2^2+
\int_\Omega Tz^2\right)\qquad(z\in H^1(\Omega)).
\label{eq:weighted-poincare}
\end{equation}
\end{lemma}

\begin{proof}
Otherwise a normalized counterexample sequence has vanishing gradient and
vanishing weighted term.  Rellich compactness and connectedness make its
strong \(L^2\) limit a nonzero constant, contradicting
\(T\geq0\) and \(\int_\Omega T>0\).
\end{proof}

For \(0\leq c,v\leq1\), denote the solutions of
\eqref{eq:elliptic-o-compact} and \eqref{eq:elliptic-a-compact} by
\(\mathcal O(c,v)\) and \(\mathcal A(c,o)\), respectively.  These maps
are the analytical replacement for fast time regularity.

\begin{proposition}[Signal resolvents]
\label{prop:resolvents}
Under (A1), (A3)--(A5), each admissible \((c,v)\) determines a unique
\(\mathcal O(c,v)\in H^1\cap L^\infty\), and each admissible \((c,o)\)
determines a unique \(\mathcal A(c,o)\in H^1\cap L^\infty\), with
\[
0\leq\mathcal O(c,v)\leq o_b,\qquad
0\leq\mathcal A(c,o)\leq\sigma_a/\mu_a.
\]
For a constant \(C\) depending on fixed parameters, \(\eps\), and
\(\ell\),
\begin{align}
\norm{\mathcal O(c_1,v_1)-\mathcal O(c_2,v_2)}_{H^1}
&\leq C(\norm{c_1-c_2}_2+\norm{v_1-v_2}_2),
\label{eq:O-stability}\\
\norm{\mathcal A(c_1,o_1)-\mathcal A(c_2,o_2)}_{H^1}
&\leq C(\norm{c_1-c_2}_2+\norm{o_1-o_2}_2).
\label{eq:A-stability}
\end{align}
Both maps preserve strong time measurability.
\end{proposition}

\begin{proof}
We first treat the oxygen resolvent. Extend \(h_K(s)=s/(K+s)\) from \([0,o_b]\) to a globally Lipschitz,
increasing function by using slope \(1/K\) below zero and the tangent
slope \(d_K=K/(K+o_b)^2\) above \(o_b\). Its slope is bounded below by \(d_K\), in the sense that
\[
(\widehat h_K(r)-\widehat h_K(s))(r-s)
\geq d_K|r-s|^2,
\qquad r,s\in\mathbb R.
\label{eq:extended-monotonicity}
\]

For fixed \(0\leq c,v\leq1\), define
\(\mathcal L^o_{c,v}:H^1(\Omega)\to(H^1(\Omega))^*\) by
\begin{equation*}
\langle\mathcal L^o_{c,v}u,\varphi\rangle
=D_o\int_\Omega\nabla u\cdot\nabla\varphi
+\gamma_T\int_\Omega T\widehat h_{K_T}(u)\varphi
+\gamma_c\int_\Omega c\,\widehat h_{K_c}(u)\varphi
+\lambda\int_\Omega\Ffun[v]\,u\varphi .
\end{equation*}
The extended oxygen problem is $\mathcal L^o_{c,v}o=\lambda o_b\Ffun[v]$ in $(H^1(\Omega))^*$.

The operator \(\mathcal L^o_{c,v}\) is well defined and bounded on bounded
subsets of \(H^1\), since \(c,\Ffun[v]\in L^\infty(\Omega)\), while the
functions \(\widehat h_K\) are globally Lipschitz.  It is also hemicontinuous:
for \(u,\varphi,\psi\in H^1(\Omega)\), the map $s\longmapsto \langle\mathcal L^o_{c,v}(u+s\varphi),\psi\rangle$ is continuous by the continuity of the scalar functions
\(\widehat h_{K_T}\) and \(\widehat h_{K_c}\) and dominated convergence.

For \(u,w\in H^1(\Omega)\), put \(z=u-w\). By the weighted Poincar\'e inequality
\eqref{eq:weighted-poincare}:
\begin{equation*}
\langle\mathcal L^o_{c,v}u-\mathcal L^o_{c,v}w,z\rangle
\geq D_o\norm{\nabla z}_2^2
+\gamma_Td_{K_T}\int_\Omega T|z|^2\geq
\alpha_o\norm{z}_{H^1}^2,
\end{equation*}
where $\alpha_o
=
\frac12\min\{D_o,\gamma_Td_{K_T}\}
\min\{1,C_{\mathrm{wp}}^{-1}\}>0$.
Thus \(\mathcal L^o_{c,v}\) is strongly monotone and coercive.  The
Browder--Minty theorem therefore yields a unique extended solution
\(o\in H^1(\Omega)\).

It remains to show that the extension is inactive. Testing the extended equation by \(-o^-\) gives
\[
D_o\norm{\nabla o^-}_2^2
+\gamma_T\int_\Omega T\widehat h_{K_T}(o)(-o^-)
+\gamma_c\int_\Omega c\widehat h_{K_c}(o)(-o^-)
+\lambda\int_\Omega\Ffun[v](o^-)^2
=
-\lambda o_b\int_\Omega\Ffun[v]o^-.
\]
Then \(o^-=0\), and therefore $o\geq0$.
Next test by \((o-o_b)^+\):
\begin{align*}
&D_o\norm{\nabla (o-o_b)^+}_2^2
+\gamma_T\int_\Omega
T\widehat h_{K_T}(o)(o-o_b)^+\\
&\quad
+\gamma_c\int_\Omega
c\widehat h_{K_c}(o)(o-o_b)^+
+\lambda\int_\Omega
\Ffun[v](o-o_b)^+{}^2
=0.
\end{align*}
$(o-o_b)^+=0$, and thus $0\leq o\leq o_b$.
Consequently \(\widehat h_K(o)=h_K(o)\), and the extended solution solves the
original oxygen equation.  This proves existence, uniqueness, and the
asserted invariant interval without any positive lower bound on
\(\Ffun[v]\).

Similarly, the VEGF bilinear form is coercive with constant
\(\min\{D_a,\mu_a\}\). Lax--Milgram Theorem and the analogous two truncation
tests ($-a^-$ for $a\geq0$ and $(a-\sigma_a/\mu_a)^+$ for $a\leq\dfrac{\sigma_a}{\mu_a}$) prove its assertions. 

We now prove \eqref{eq:O-stability}.  Write
\[
o_i=\mathcal O(c_i,v_i),\qquad
F_i=\Ffun[v_i],\qquad
\delta o=o_1-o_2.
\]
Subtracting the two oxygen equations and retaining on the left the
monotone terms associated with \((c_1,F_1)\), we obtain
\begin{align*}
&D_o\int_\Omega\nabla\delta o\cdot\nabla\varphi
+\gamma_T\int_\Omega
T\bigl[\widehat h_{K_T}(o_1)-\widehat h_{K_T}(o_2)\bigr]\varphi\\
&\quad
+\gamma_c\int_\Omega
c_1\bigl[\widehat h_{K_c}(o_1)-\widehat h_{K_c}(o_2)\bigr]\varphi
+\lambda\int_\Omega F_1\delta o\,\varphi\\
&=
\lambda o_b\int_\Omega(F_1-F_2)\varphi
-\gamma_c\int_\Omega(c_1-c_2)\widehat h_{K_c}(o_2)\varphi
-\lambda\int_\Omega(F_1-F_2)o_2\varphi .
\end{align*}
The two terms involving \(F_1-F_2\) are combined before estimating:
\begin{equation*}
\lambda o_b(F_1-F_2)
-\lambda(F_1-F_2)o_2
=
\lambda(F_1-F_2)(o_b-o_2).
\end{equation*}
Thus the right-hand side becomes
\[
\lambda(F_1-F_2)(o_b-o_2)
-\gamma_c(c_1-c_2)\widehat h_{K_c}(o_2).
\]
Testing by \(\varphi=\delta o\), using
\(0\leq o_2\leq o_b\) and
\(0\leq\widehat h_{K_c}(o_2)\leq1\), gives
\[
\begin{aligned}
\alpha_o\norm{\delta o}_{H^1}^2
&\leq
\left[
\gamma_c\norm{c_1-c_2}_2
+\lambda o_b\norm{F_1-F_2}_2
\right]\norm{\delta o}_2\\
&\leq
\left[
\gamma_c\norm{c_1-c_2}_2
+\lambda o_bL_{\Ffun}\norm{v_1-v_2}_2
\right]\norm{\delta o}_{H^1}.
\end{aligned}
\]
Hence
\[
\norm{\mathcal O(c_1,v_1)-\mathcal O(c_2,v_2)}_{H^1}
\leq
\frac{1}{\alpha_o}
\left[
\gamma_c\norm{c_1-c_2}_2
+\lambda o_bL_{\Ffun}\norm{v_1-v_2}_2
\right],
\]
which proves \eqref{eq:O-stability}. Similar procedures with 
\[
\operatorname{Lip}(H_H)\leq K_H^{-1},
\qquad
0\leq a_2\leq\frac{\sigma_a}{\mu_a},
\]
and the coercivity estimate above yields
\[
\norm{\mathcal A(c_1,o_1)-\mathcal A(c_2,o_2)}_{H^1}
\leq
C\left(
\norm{c_1-c_2}_2+\norm{o_1-o_2}_2
\right),
\]
which proves \eqref{eq:A-stability}.

The stability estimates show that $(c,v)\mapsto\mathcal O(c,v)$ and $(c,o)\mapsto\mathcal A(c,o)$
are continuous maps between separable Banach spaces.  Hence they preserve
strong measurability under composition with strongly measurable
time-dependent inputs.
\end{proof}

A quasi-static weak solution has the slow-variable regularity in
\eqref{eq:weak-reg-diffusive}--\eqref{eq:weak-reg-ode} (without the fast
time derivatives), satisfies \eqref{eq:weak-c-compact}--\eqref{eq:weak-v-compact}, and satisfies
\eqref{eq:elliptic-o-compact}--\eqref{eq:elliptic-a-compact} weakly at
almost every time.  Only the slow fields carry initial conditions.

\begin{theorem}[Fast-signal limit]
\label{thm:fast}
Assume (A1)--(A7), fix \(\eps,\ell>0\), and let
\((\tau_o^n,\tau_a^n)\to(0,0)\).  Any sequence of global weak solutions
with common admissible initial data has a subsequence that converges on every
finite interval to a global quasi-static weak solution.  Precisely,
\begin{align}
c_n,p_n&\rightharpoonup c,p&&\text{in }L^2(0,T_f;H^1),
&c_n,p_n&\to c,p&&\text{in }L^2(Q_T),
\label{eq:fast-cp}\\
f_n&\to f&&\text{in }C([0,T_f];L^2),
&v_n&\to v&&\text{in }C([0,T_f];L^2),
\label{eq:fast-fv}\\
o_n&\to o&&\text{in }L^2(0,T_f;H^1),
&a_n&\to a&&\text{in }L^2(0,T_f;H^1),
\label{eq:fast-ao}\\
\pi_n&\to\pi&&\text{in }L^\infty(0,T_f;H^1),
&(q[v_n],\Ffun[v_n])&\to(q[v],\Ffun[v])
&&\text{in }L^\infty(0,T_f;L^2)^2.
\label{eq:fast-perfusion}
\end{align}
The strong convergences in \eqref{eq:fast-cp}--\eqref{eq:fast-ao} also hold
almost everywhere after extraction; bounded components converge weak-star
in \(L^\infty\), and \(f_n,v_n\) retain their weak-star derivative
bounds.  Moreover
\begin{equation}
\tau_o^n(o_n)_t\rightharpoonup0,\qquad
\tau_a^n(a_n)_t\rightharpoonup0
\quad\text{in }L^2(0,T_f;(H^1)^*),
\label{eq:fast-residuals}
\end{equation}
and \(o=\mathcal O(c,v)\), \(a=\mathcal A(c,o)\) almost everywhere.
The theorem does not assert uniqueness for the full limiting system.
\end{theorem}

The proof of Theorem~\ref{thm:fast} is provided in Appendix~\ref{app:fast-details}. A quantitative \(O(\sqrt{\tau_o}+\sqrt{\tau_a})\) rate would require a
unique regular solution family and well-prepared data (or an explicit
initial-layer term).  Since those hypotheses are not proved here, no rate
theorem is claimed.

\section{Reduced invasion fronts}
\label{sec:waves}
This section concerns a local scalar reduction, not a traveling wave of the
complete pressure-coupled system.  Freeze the one-dimensional taxis gradients
and oxygen in a patch, neglect transverse variation, and assume the maturation
loss is absent there (for example \(v_{\rm fr}=1\)).  With
\begin{equation}
B=\chi_a\partial_xS(a)|_{\rm fr}+\chi_f\partial_xf|_{\rm fr},
\qquad D=D_c,\qquad R=r_cH_P(o_{\rm fr})>0,
\label{eq:frozen-coefficients}
\end{equation}
the endothelial equation reduces to
\begin{equation}
c_t=Dc_{xx}-B\partial_x[c(1-c)]+Rc(1-c).
\label{eq:fisher-burgers}
\end{equation}
Here \(B>0\) is forward taxis. 
A decreasing front
\(c(x,t)=\phi(\xi)\), \(\xi=x-st\), with \(\phi(-\infty)=1\) and
\(\phi(+\infty)=0\), obeys
\begin{equation}
D\phi''+[s-B(1-2\phi)]\phi'+R\phi(1-\phi)=0.
\label{eq:wave-profile}
\end{equation}
As \(\xi\to+\infty\), the decreasing front satisfies
\(\phi(\xi)\to0\).  Linearizing \eqref{eq:wave-profile} at the
equilibrium \(\phi=0\) gives $D\phi''+(s-B)\phi'+R\phi=0$.
Seeking a decaying exponential mode
\(\phi(\xi)\sim A e^{-\rho\xi}\), \(\rho>0\), yields the characteristic
equation $D\rho^2-(s-B)\rho+R=0$.
Hence the leading-edge analysis imposes the necessary condition $s\ge B+2\sqrt{DR}$.
At the double-root case \(s=B+2\sqrt{DR}\), the linearized modes have
the form \((A_1+A_2\xi)e^{-\sqrt{R/D}\xi}\); thus the pure exponential
ansatz is understood as the generic leading-order mode.
To remove the constant drift \(B\), set $y=x-Bt$, $u(y,t)=c(x,t)$, $\zeta=\sqrt{\frac RD}\,y$, $\vartheta=Rt$, we obtain
\begin{equation}
u_\vartheta+\mu uu_\zeta=u_{\zeta\zeta}+u(1-u),\quad \mu=-\frac{2B}{\sqrt{DR}}, \quad
s=B+\sqrt{DR}\,\sigma.
\label{eq:scaled-burgers}
\end{equation}
The equation for $s$ in the \((\zeta,\vartheta)\)-variables comes from the fact that the speed of a front of physical speed \(s\) is $\sigma=\dfrac{s-B}{\sqrt{DR}}$.

\begin{lemma}[Exact Fisher--Burgers threshold]
\label{lem:burgers}
The dimensionless equation in \eqref{eq:scaled-burgers} has a monotone
front of speed \(\sigma\) if and only if \(\sigma\geq\sigma_*(\mu)\),
where
\begin{equation}
\sigma_*(\mu)=
\begin{cases}2,&\mu\leq2,\\
\mu/2+2/\mu,&\mu>2.
\end{cases}
\label{eq:dimensionless-speed}
\end{equation}
For \(\mu>2\), the threshold speed is attained by the explicit profile
\[
U(z)=\frac{1}{1+e^{\mu z/2}},
\qquad
z=\zeta-\sigma_*(\mu)\vartheta.
\]
\end{lemma}

\begin{proof}
We write \eqref{eq:scaled-burgers} in the form
\[
u_\vartheta+h_\mu(u)u_\zeta=u_{\zeta\zeta}+f(u),
\qquad
h_\mu(u)=\mu u,\qquad
f(u)=u(1-u).
\]
Then \(h_\mu\in C^1([0,1])\), \(f\in C^\infty([0,1])\), and
\[
f(0)=f(1)=0,\qquad
f(u)>0\quad(0<u<1),\qquad
f'(0)=1>0>f'(1)=-1.
\]
Thus the standard threshold result for monotone travelling fronts of
reaction--diffusion equations with nonlinear convection
\cite{malaguti2002travelling} applies.  In particular, there is a critical
speed \(\sigma_*(\mu)\) such that monotone decreasing fronts
connecting \(1\) to \(0\) occur precisely at and above this threshold.
It remains to determine \(\sigma_*(\mu)\) explicitly.

For a travelling front
\[
u(\zeta,\vartheta)=U(z),\qquad z=\zeta-\sigma\vartheta,
\]
with \(U(-\infty)=1\), \(U(+\infty)=0\), and \(U'<0\), the profile satisfies
\begin{equation}
U''+(\sigma-\mu U)U'+U(1-U)=0.
\label{eq:burgers-profile-proof}
\end{equation}
At the leading edge \(U=0\), linearization gives $U''+\sigma U'+U=0$.
A decaying mode \(U(z)\sim e^{-\rho z}\), \(\rho>0\), therefore satisfies $\rho^2-\sigma\rho+1=0$. Hence every monotone front necessarily satisfies $\sigma\ge2$.

We next use the variational characterization of the minimal speed for
reaction--convection--diffusion fronts established in
\cite{benguria2004minimal}.  For \(\mu\ne0\), write
\[
\widetilde\mu=|\mu|,
\qquad
\varphi_\mu(u)=\operatorname{sgn}(\mu)\,u,
\]
so that
\[
\widetilde\mu\,\varphi_\mu(u)=\mu u,
\qquad
\varphi_\mu\in C^1([0,1]),
\qquad
\varphi_\mu(0)=0.
\]
For \(\mu=0\), the corresponding formula reduces to the classical
Fisher--KPP variational principle.  Thus, with the admissible class
\[
\mathcal G
=
\left\{
g\in C^1((0,1)):
g>0,\ g'<0,\ g(1)=0,
\ \text{and all terms below are integrable}
\right\},
\]
the minimal speed is characterized by
\begin{equation}
\sigma_*(\mu)
=
\sup_{g\in\mathcal G}
\mathcal E_\mu(g),
\qquad
\mathcal E_\mu(g)
=
\frac{
\displaystyle
\int_0^1
\left[
2\sqrt{u(1-u)g(-g')}
+\mu ug
\right]\dd u
}{
\displaystyle
\int_0^1g\dd u
}.
\label{eq:burgers-variational-proof}
\end{equation}

We first derive an upper bound.  Fix \(a>0\).  By the arithmetic--geometric
mean inequality,
\[
\begin{aligned}
\int_0^1
2\sqrt{u(1-u)g(-g')}\,\dd u
&\le
a\int_0^1g\dd u
+\frac1a
\int_0^1u(1-u)(-g')\dd u.
\end{aligned}
\]
Integration by parts gives
\[
\int_0^1u(1-u)(-g')\dd u =
-\bigl[u(1-u)g(u)\bigr]_{0}^{1}
+\int_0^1(1-2u)g(u)\dd u=
\int_0^1(1-2u)g(u)\dd u.
\]
Consequently,
\[
\mathcal E_\mu(g)\le
\frac{
\displaystyle
\int_0^1
\left[
a+\frac{1-2u}{a}+\mu u
\right]g(u)\dd u
}{
\displaystyle
\int_0^1g(u)\dd u
}\le
\sup_{0\le u\le1}
\left[
a+\frac{1-2u}{a}+\mu u
\right].
\]
Since \(a>0\) is arbitrary,
\begin{equation}
\sigma_*(\mu)
\le
\inf_{a>0}
\sup_{0\le u\le1}
\left[
a+\frac{1-2u}{a}+\mu u
\right].
\label{eq:burgers-upper-proof}
\end{equation}

Suppose first that \(\mu\le2\).  Choosing \(a=1\) in
\eqref{eq:burgers-upper-proof} gives
\[
\sup_{0\le u\le1}
[1+(1-2u)+\mu u]
=
\sup_{0\le u\le1}
[2+(\mu-2)u]
=2.
\]
Hence $\sigma_*(\mu)\le2$.
Combined with the necessary leading-edge condition
\eqref{eq:burgers-profile-proof}, this yields $\sigma_*(\mu)=2$, $\mu\le2.$

Now suppose that \(\mu>2\).  Choosing $a=\frac2\mu$ in \eqref{eq:burgers-upper-proof} gives
\begin{equation}
\sigma_*(\mu)
\le
\frac{\mu}{2}+\frac2\mu.
\label{eq:burgers-upper-pushed}
\end{equation}

To obtain the matching lower bound, choose
\[
g_\alpha(u)
=
\left(\frac{1-u}{u}\right)^\alpha,
\qquad
0<\alpha<1.
\]
Then \(g_\alpha>0\), \(g_\alpha'<0\), \(g_\alpha(1)=0\), and $g_\alpha(u)\sim u^{-\alpha}$ as $u\downarrow0$, so \(\int_0^1g_\alpha\,\dd u<\infty\).  Moreover, $2\sqrt{u(1-u)g_\alpha(-g_\alpha')}
=
2\sqrt{\alpha}\,g_\alpha$.
It follows that
\[
\mathcal E_\mu(g_\alpha)
=
2\sqrt{\alpha}
+
\mu
\frac{\displaystyle\int_0^1u g_\alpha(u)\dd u}
{\displaystyle\int_0^1g_\alpha(u)\dd u}.
\]
Since $g_\alpha(u)=u^{-\alpha}(1-u)^\alpha$,
the Beta-function identity gives $\int_0^1g_\alpha(u)\dd u
=
B(1-\alpha,1+\alpha)$,
and $\int_0^1u g_\alpha(u)\dd u
=
B(2-\alpha,1+\alpha)
=
\frac{1-\alpha}{2}
B(1-\alpha,1+\alpha)$.
Thus $\mathcal E_\mu(g_\alpha)
=
2\sqrt{\alpha}
+\frac{\mu}{2}(1-\alpha)$.
Taking $\alpha=\frac4{\mu^2}\in(0,1)$ gives $\mathcal E_\mu(g_{4/\mu^2}) =
\frac{\mu}{2}+\frac2\mu$.
Hence, by \eqref{eq:burgers-variational-proof} and \eqref{eq:burgers-upper-pushed}, $\sigma_*(\mu)
=
\frac{\mu}{2}+\frac2\mu$ for $\mu>2$.

It remains to verify the explicit minimal profile. Set $a=\frac{\mu}{2}$ and $U(z)=\dfrac{1}{1+e^{az}}$.
Then $U'=-aU(1-U)$ and $U''=a^2U(1-U)(1-2U)$.
Substitution into \eqref{eq:burgers-profile-proof} gives
\[
\begin{aligned}
U''+(\sigma-\mu U)U'+U(1-U)
&=
U(1-U)
\bigl[
a^2(1-2U)-a(\sigma-\mu U)+1
\bigr].
\end{aligned}
\]
Because \(a=\mu/2\), to let the bracketed term vanish we set 
$\sigma=a+\frac1a
=
\frac{\mu}{2}+\frac2\mu$, and thus
$U(z)=\frac{1}{1+e^{\mu z/2}}$
is a monotone decreasing front at the threshold speed for every
\(\mu>2\).

Therefore
\[
\sigma_*(\mu)=
\begin{cases}
2,&\mu\le2,\\[1mm]
\dfrac{\mu}{2}+\dfrac2\mu,&\mu>2,
\end{cases}
\]
and the asserted explicit minimal profile holds in the pushed regime.
The travelling-front threshold theorem then yields existence of monotone
fronts for every \(\sigma\ge\sigma_*(\mu)\) and excludes monotone fronts
below the threshold.
\end{proof}

\begin{figure}[htbp]
\centering
\includegraphics[width=\linewidth]{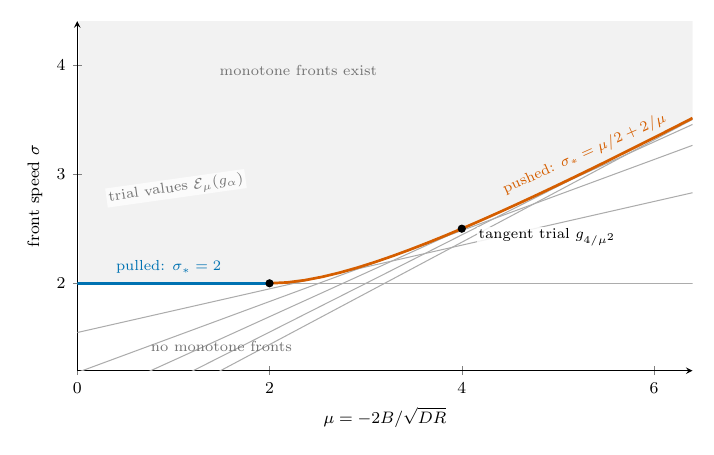}
\caption{Lemma~\ref{lem:burgers} as an envelope construction.  Thin gray
lines are the trial values
\(\mathcal E_\mu(g_\alpha)=2\sqrt\alpha+\tfrac\mu2(1-\alpha)\) for a
family of admissible \(g_\alpha\); their upper envelope (thick) is the
exact threshold \(\sigma_*(\mu)\) in \eqref{eq:dimensionless-speed}.
Monotone fronts
exist exactly in the shaded region \(\sigma\geq\sigma_*(\mu)\); the
pushed branch is tangent to the trial line with \(\alpha=4/\mu^2\).}
\label{fig:variational}
\end{figure}

Figure~\ref{fig:variational} shows the variational mechanism behind the
two branches: the threshold is the upper envelope of the trial values
\(\mathcal E_\mu(g_\alpha)\), which are affine in \(\mu\); for
\(\mu\leq2\) the envelope is pinned at the leading-edge value \(2\),
while for \(\mu>2\) the tangent trial function \(g_{4/\mu^2}\) realizes
the strictly larger pushed value.

The resulting physical theorem exposes a nonlinear branch that
leading-edge linearization alone misses.

\begin{theorem}[Reduced invasion-front speed]
\label{thm:wave}
Under the local freezing assumptions preceding \eqref{eq:frozen-coefficients},
equation \eqref{eq:fisher-burgers} admits a monotone invasion front of speed
\(s\) if and only if \(s\geq s_*(B,R)\), where
\begin{equation}
s_*(B,R)=
\begin{cases}
B+2\sqrt{DR},&B\geq-\sqrt{DR},\\[1mm]
-DR/B,&B<-\sqrt{DR}.
\end{cases}
\label{eq:physical-speed}
\end{equation}
\end{theorem}

\begin{proof}
By Lemma~\ref{lem:burgers}, the scaled equation admits a monotone
decreasing front of speed \(\sigma\) if and only if $\sigma\ge \sigma_*(\mu)$ with $\sigma_*$ given by \eqref{eq:dimensionless-speed}.
The parameters in the scaled equation are related to the physical
parameters by $\mu=-\frac{2B}{\sqrt{DR}}$ and $s=B+\sqrt{DR}\,\sigma$.
Since \(D>0\) and \(R>0\), the map
\[
\sigma\longmapsto s=B+\sqrt{DR}\,\sigma
\]
is a strictly increasing bijection, with inverse $\sigma=\dfrac{s-B}{\sqrt{DR}}$.
Hence the physical minimal speed is
\[
s_*(B,R)
=
B+\sqrt{DR}\,\sigma_*(\mu),
\qquad
\mu=-\frac{2B}{\sqrt{DR}}.
\]

If \(\mu\le2\), then
\[
-\frac{2B}{\sqrt{DR}}\le2
\quad\Longleftrightarrow\quad
B\ge-\sqrt{DR},
\]
and therefore $s_*(B,R)
=
B+2\sqrt{DR}$.

If \(\mu>2\), then
\[
-\frac{2B}{\sqrt{DR}}>2
\quad\Longleftrightarrow\quad
B<-\sqrt{DR}.
\]
In this case, $s_*(B,R)=
B+\sqrt{DR}
\left(
\frac{\mu}{2}+\frac{2}{\mu}
\right)=
-\dfrac{DR}{B}$.
Thus
\[
s_*(B,R)=
\begin{cases}
B+2\sqrt{DR},&B\ge-\sqrt{DR},\\[1mm]
-DR/B,&B<-\sqrt{DR}.
\end{cases}
\]

Finally, \(s_*(B,R)\) is continuous at the branch point \(B=-\sqrt{DR}\).  Since the
coordinate transformation between \((x,t)\) and \((\zeta,\vartheta)\) is
invertible, the existence and nonexistence statements for monotone
fronts transfer equivalently between the scaled and physical equations.
Hence \eqref{eq:physical-speed} holds.
\end{proof}

For pre- and post-perfusion frozen states, put
\(R_\pm=r_cH_P(o_\pm)>0\).  If both are on the linear branch, subtraction
of \eqref{eq:physical-speed} gives
\begin{equation*}
s_+-s_-=(B_+-B_-)+2\sqrt{Dr_c}
[\sqrt{H_P(o_+)}-\sqrt{H_P(o_-)}],
\end{equation*}
and hence acceleration exactly when
\begin{equation*}
B_+-B_->-2\sqrt{Dr_c}
[\sqrt{H_P(o_+)}-\sqrt{H_P(o_-)}].
\end{equation*}
If either state is on the adverse-taxis branch, the two cases in
\eqref{eq:physical-speed} must be applied separately.  This is a local
\emph{speed-transition criterion}, not a bifurcation or a clinical
prediction; the full model contains neither a translating nonlocal state nor
a treatment variable.

\section{Numerical simulations}
\label{sec:numerical}

\subsection{Fully specified numerical baseline}
\label{subsec:numerical-scheme}
On a Cartesian face between cells \(i\) and \(i+1\), the discrete
endothelial flux in one coordinate direction is
\begin{equation}
\mathcal Q^c_{i+1/2}=D_c\frac{c_{i+1}-c_i}{h}
-\overline M_{i+1/2}\left[
\chi_a\frac{S(a_{i+1})-S(a_i)}{h}
+\chi_f\frac{f_{i+1}-f_i}{h}\right],
\label{eq:discrete-c-flux}
\end{equation}
where \(\overline M_{i+1/2}=[M(c_i)+M(c_{i+1})]/2\).  The cell residual
is the difference between the complete face fluxes.  This formula is
conservative
and second-order consistent on the smooth manufactured test. Pressure uses
harmonic interior transmissibilities
\(K_{i+1/2}=2K_iK_{i+1}/(K_i+K_{i+1})\) and half-cell distances on the
Dirichlet faces.  At cell centers \(x_i\),
\begin{equation}
(\mathcal J_{\ell,h}v)_i=
\frac{\sum_jJ_\ell(x_i-x_j)v_j\abs{C_j}}
{\sum_jJ_\ell(x_i-x_j)\abs{C_j}}.
\label{eq:discrete-kernel}
\end{equation}
Identical numerator and denominator stencils preserve constants at physical
boundaries.
One time step performs the following ordered operations: compute
\(\pi^n,q^n,\Ffun[v^n]\); update \(p\) by backward Euler for diffusion and
decay; solve the nonlinear backward-Euler oxygen equation and then the VEGF
equation; update \(c\)
with backward-Euler diffusion and the conservative explicit taxis flux
\eqref{eq:discrete-c-flux}, in which the taxis gradients are evaluated
with the just-updated \(a^{n+1}\) and the lagged \(f^{n}\); and
apply the rational one-step formulas for the ECM and vessel reactions,
\[
f^{n+1}=\frac{f^n+\Delta t\,\rho_f}
{1+\Delta t(\rho_f+\kappa_fp^{n+1})},\qquad
v^{n+1}=\frac{v^n+\Delta t\,\eta c^{n+1}}
{1+\Delta t\{\eta c^{n+1}+\delta_v[1-G(q^n)]\}},
\]
which preserve the invariant intervals unconditionally.
Pressure is therefore evaluated at the single explicit stage, not lagged
inside a multistage method. 

\subsection{Auxiliary numerical verification}
\label{subsec:numerics}
All computations are nondimensional and serve as verification and mechanism
tests, not as parameter identification.  We use conservative cell-centered
finite volumes on rectangular meshes.  The complete endothelial flux is
assembled facewise; pressure uses harmonic face conductivity and the same
boundary-normalized kernel as \eqref{eq:average}.  Backward Euler treats
diffusion and linear relaxation, while the remaining terms use an IMEX Euler
split; pressure is recomputed at the only stage. Higher-order IMEX schemes require full-precision coefficients
and separate positivity analysis for the taxis--pressure coupling
\cite{de2026imex}.
Table~\ref{tab:verification} gives the verified results.  Spatial
and temporal refinement were performed separately, and algebraic errors were
made negligible relative to discretization errors.  The kernel entry is
exactly zero, not merely small: the test field is constant, so the discrete
numerator and denominator in \eqref{eq:discrete-kernel} consist of
identical terms summed in the identical order, and the quotient equals one.

\begin{table}[htbp]
\centering
\caption{Numerical verification. Reported rates are from the last refinement.}
\label{tab:verification}
\begin{tabularx}{\textwidth}{@{} >{\raggedright\arraybackslash}X l >{\raggedright\arraybackslash}X @{}}
\toprule
Test & Result & Acceptance condition\\
\midrule
Constant-\(v\) pressure & \(L^2\) error \(2.34\times10^{-14}\) & \(<10^{-8}\)\\
Normalized kernel, including corners & maximum error \(0\) & \(<5\times10^{-14}\)\\
Manufactured taxis-flux spatial test & rate \(1.855\) & \(\abs{r-2}\leq0.2\)\\
Backward-Euler temporal test & rate \(1.096\) & \(\abs{r-1}\leq0.2\)\\
Pressure-solver maximum residual & \(1.3\times10^{-15}\) & \(<10^{-12}\)\\
\bottomrule
\end{tabularx}
\end{table}

\begin{figure}[htbp]
\centering
\includegraphics[width=0.48\textwidth]{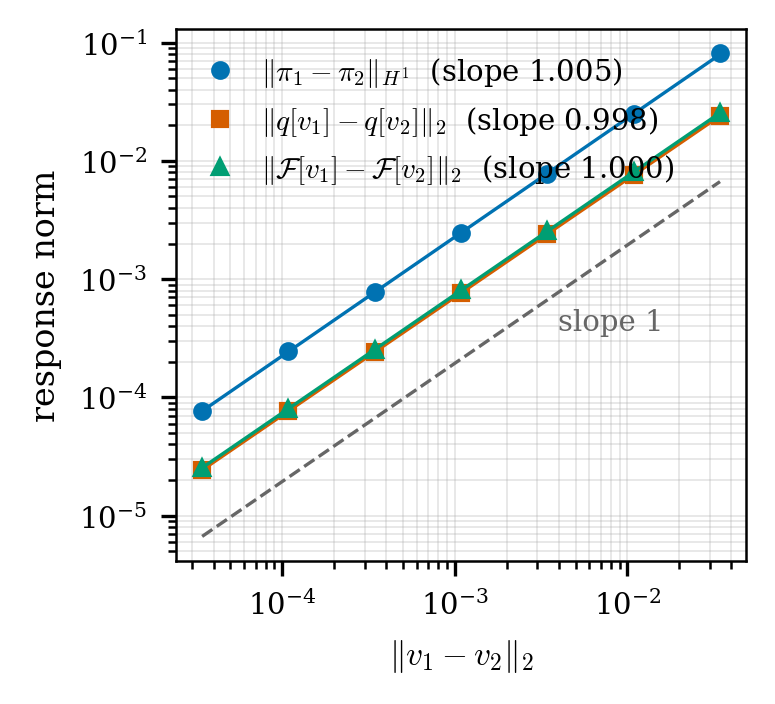}
\hfill
\includegraphics[width=0.48\textwidth]{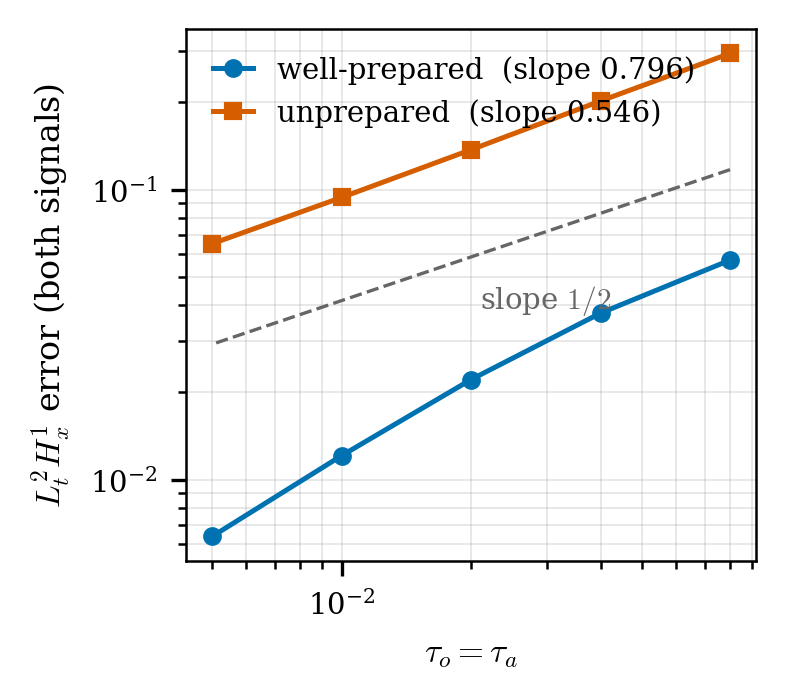}
\caption{Pressure--perfusion stability (left) and dynamic-to-quasi-static
signal error for prepared and unprepared data (right).}
\label{fig:pressure-fast}
\end{figure}

For the pressure test, vessel perturbations spanning three orders of
magnitude produce fitted log--log slopes 1.005, 0.998, and 1.000 for
\(\pi,q,\Ffun\), respectively.  This result is consistent with
Theorem~\ref{thm:perfusion}.  
In the fast-signal test, smooth time-dependent \(c,v\) are prescribed and pressure is recomputed.  With
\(\tau_o=\tau_a\in\{0.08,0.04,0.02,0.01,0.005\}\), the combined
\(L^2_tH^1_x\) error has fitted slope 0.796 for well-prepared data and
0.546 for unprepared data (Fig.~\ref{fig:pressure-fast}).  The unprepared
slope is pinned near \(\tfrac12\) by the initial layer; the prepared
curve steepens toward first order as \(\tau\) decreases but is still
pre-asymptotic at the largest \(\tau\) shown.  

\begin{figure}[htbp]
\centering
\includegraphics[width=0.76\textwidth]{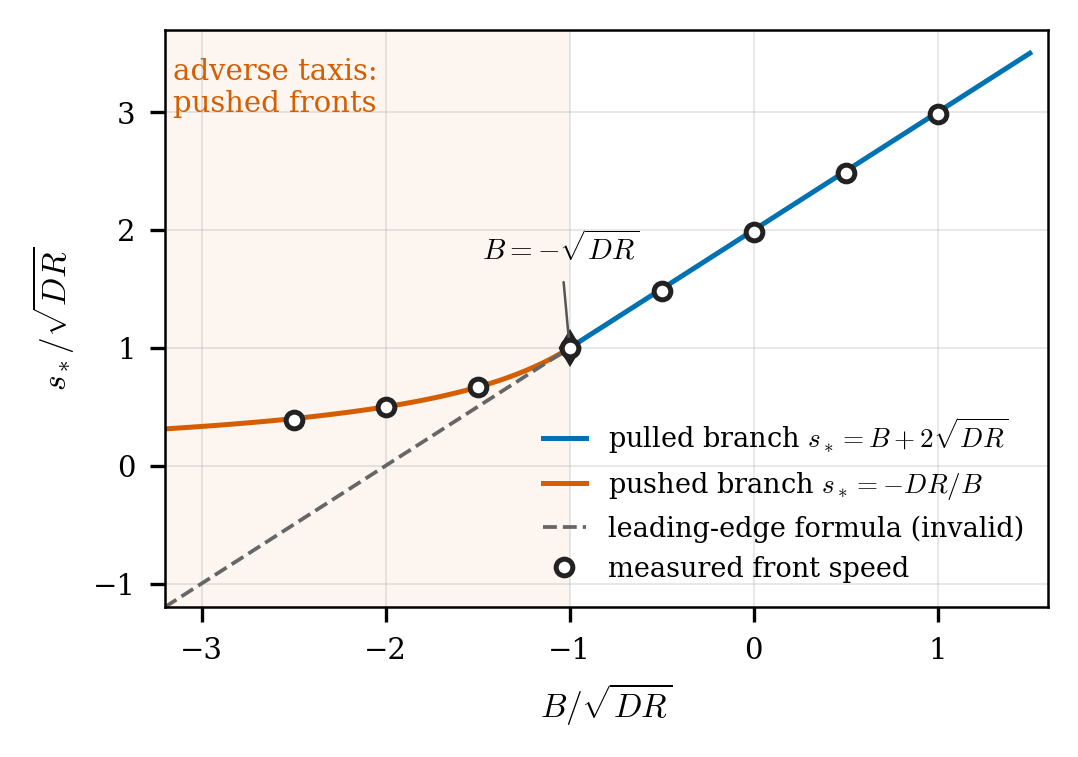}
\caption{Measured invasion-front speeds (open markers) against the exact
piecewise threshold of Theorem~\ref{thm:wave}: the pulled branch
\(s_*=B+2\sqrt{DR}\), the pushed branch \(s_*=-DR/B\) in the
adverse-taxis region (shaded), and the leading-edge formula continued
past \(B=-\sqrt{DR}\), where it is invalid (dashed).}
\label{fig:waves}
\end{figure}

For \eqref{eq:fisher-burgers}, the front speed is fitted from the
\(c=0.5\) level set after discarding the transient, for taxis strengths
\(B/\sqrt{DR}\in[-2.5,1]\) covering both branches of
\eqref{eq:physical-speed}; steep pushed fronts (\(\mu>2\)) are computed
on a finer mesh from the exact minimal profile.  The measured speeds
match the exact threshold within 2.45\% at the steepest pushed front
(\(B=-2.5\sqrt{DR}\)) and within 1.02\% at every other tested value on
either branch (Fig.~\ref{fig:waves}).

\begin{figure}[htbp]
\centering
\includegraphics[width=0.88\textwidth]{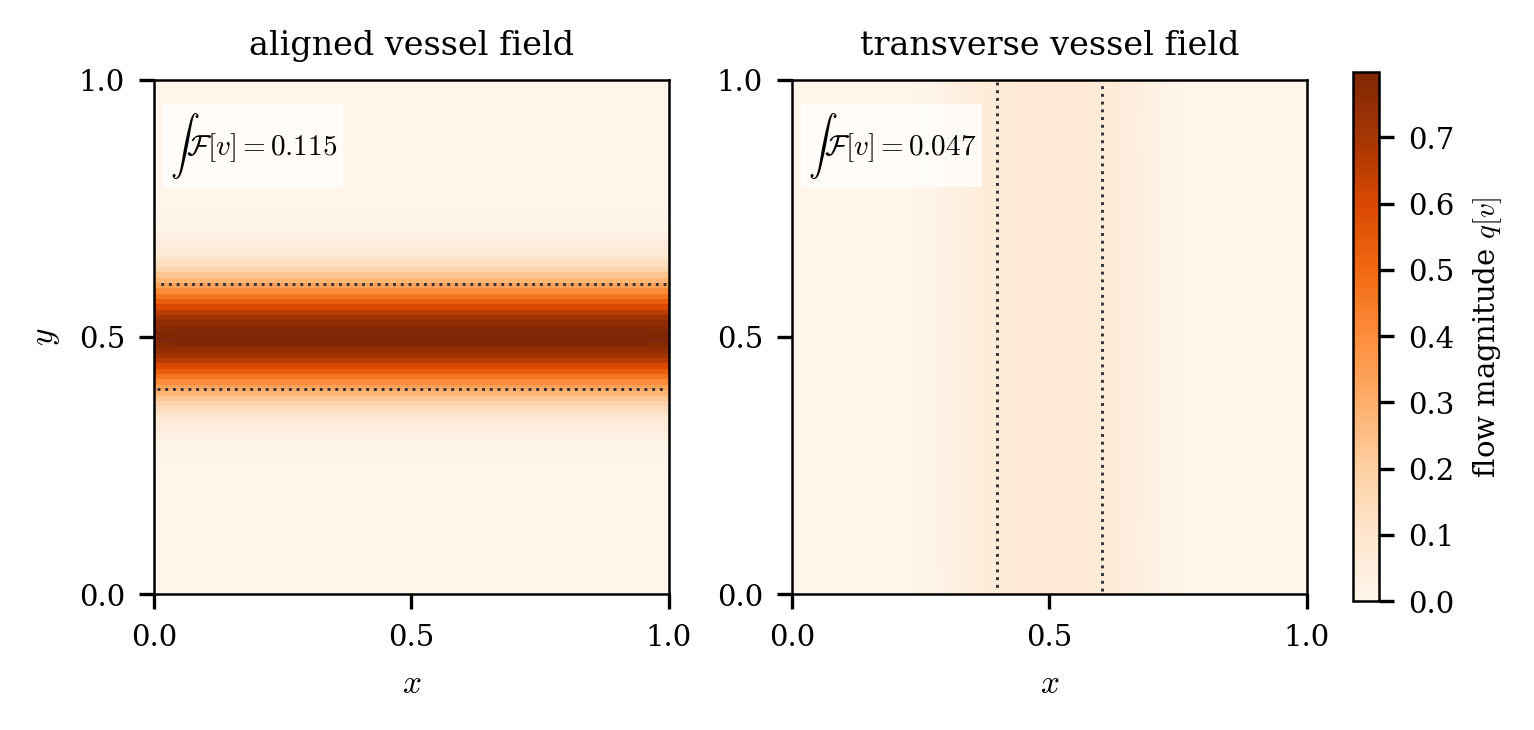}\\[1mm]
\includegraphics[width=0.88\textwidth]{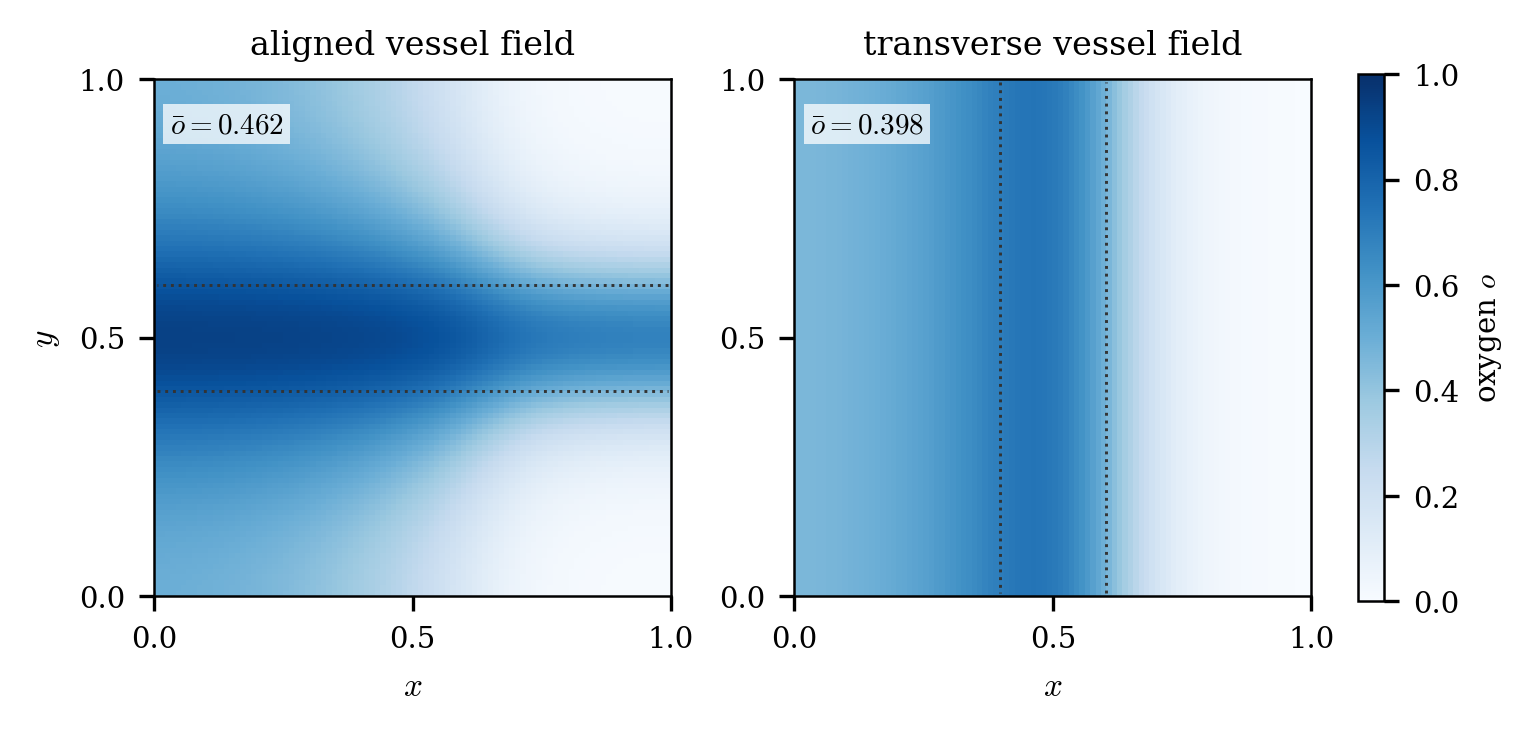}
\caption{Equal-mass flow-contrast experiment: pressure-derived flow
magnitude (top row) and oxygen (bottom row) for two vessel fields of
identical structural mass, one aligned with the arterial--venous drop and
one transverse to it (dotted contours: \(v=0.3\)).  Equal structure,
unequal function.}
\label{fig:equal-mass}
\end{figure}

Finally, two prescribed vessel fields have the same structural mass
\(\int_\Omega v=0.15066\), but different orientation relative to the
pressure drop: one spans \(\Gamma_A\) to \(\Gamma_V\), the other is a
transverse strip of equal mass touching neither Dirichlet part.  The
aligned and transverse fields have functional masses 0.11492 and
0.04675, mean oxygen 0.46239 and 0.39771, and hypoxic fractions
\(\abs{\{o<o_h\}}/\abs\Omega\) of 0.5191 and 0.5729
(Fig.~\ref{fig:equal-mass}).  The contrast persists at one coarser and
one finer mesh level. 
These controlled comparisons show that structural mass need not determine the regularized flow-functional proxy; they do not establish exact
network topology.  A \((q_*,\delta_v)\) sensitivity sweep over nine
parameter pairs additionally shows that the simulations retain 4.0--4.9
times as much vessel density in the high-flow half as in the low-flow half,
with the contrast increasing in \(\delta_v\) --- qualitatively consistent with
flow-sensitive remodeling \cite{floryan2025long}.

\section{Discussion and limitations}
\label{sec:discussion}
The modeling contribution is the explicit separation between motile
endothelium \(c\), lumenized structure \(v\), and flow-functional structure
\(\Ffun[v]\).  This allows oxygenation and regression to respond to a
globally solved pressure field while retaining a continuum state space.
Mathematically, the normalized kernel makes the coefficient map Lipschitz,
the positive background conductivity ensures uniform ellipticity, and the
vascular stability lemma compensates for the lack of spatial compactness.
The fast-limit proof then uses the coupled slow trajectories rather than a
frozen-coefficient relaxation.
The price of this closure is equally important.  \(\Ffun\) is a regularized
flow proxy, not a reconstruction of arterial--venous graph connectivity.
The pressure drop is macroscopic; radii, red-cell rheology, intravascular
advection, and explicit anastomoses are absent.  The tumor does not grow,
and tissue mechanics and treatment are not modeled.  Moreover, the constants
become singular as \(\eps\) or \(\ell\) vanishes, so neither limit is
covered.  These restrictions prevent interpreting the simulations as
patient-specific predictions.
The analysis suggests three sharply defined extensions: a singular
conductivity limit with possible disconnected perfused regions, weak--strong
uniqueness within a verifiable regularity class, and traveling coherent
structures for a coupled pressure system rather than its locally frozen
scalar diagnostic.  Each requires new analysis and is intentionally kept
outside the present claims.

\backmatter

\begin{appendices}

\section{Global bounded weak solution}
\label{app:approximation}
This appendix supplies the proof of Theorem~\ref{thm:global}.

\begin{proof}[Proof of Theorem~\ref{thm:global}]
Fix $T_f>0$. Write $H=L^2(\Omega)$, $V=H^1(\Omega)$, and
$Q_T=\Omega\times(0,T_f)$; $(\cdot,\cdot)$ denotes the $H$ inner
product. Constants may depend on $T_f$, the initial data, and the fixed
parameters, including $\tau_a,\tau_o,\varepsilon,\ell$, but not on the
approximation indices.

\medskip\noindent
\textit{Step 1: A fully specified approximation.}
Define
\[
\mathsf C(s)=\min\{1,\max\{0,s\}\},\qquad
\mathsf O(s)=\min\{o_b,\max\{0,s\}\},\qquad s^+=\max\{s,0\}.
\]
Use the mobility $M_*$ from (A7), and choose
\[
S_*(s)=\begin{cases}S(s),&s\geq0,\\s,&s<0.\end{cases}
\]
Then $S_*\in C^1(\mathbb R)$, $0<S_*'\leq1$, and $S_*'$ is globally
Lipschitz. All evaluations of $H_P,H_H$, and $s/(K+s)$ below have
nonnegative arguments.

Let $R_\nu=(I-\nu\Delta_N)^{-1}$, $\nu>0$, be the Neumann resolvent.
It preserves positivity and constants. Its spectral representation gives
\begin{equation}
\begin{aligned}
\|R_\nu z\|_2&\leq\|z\|_2 &&(z\in H),\\
\|\nabla R_\nu z\|_2&\leq\|\nabla z\|_2 &&(z\in V),\\
\|R_\nu z-z\|_2&\leq\sqrt{\nu}\,\|\nabla z\|_2 &&(z\in V).
\end{aligned}
\label{eq:resolvent-estimates}
\end{equation}
Also $\nabla R_\nu:H\to L^2(\Omega;\mathbb R^2)$ is bounded for
fixed $\nu$.

Let $\{e_j\}_{j\geq1}$ be an $H$-orthonormal Neumann eigenbasis,
$X_n=\operatorname{span}\{e_1,\ldots,e_n\}$, and $P_n$ the
$H$-orthogonal projection onto $X_n$. The eigenfunctions belong to
$H^2(\Omega)\subset L^\infty(\Omega)$, and
\[
\|P_n\varphi\|_V\leq\|\varphi\|_V,\qquad
P_n\varphi\to\varphi\quad\text{in }V\quad(\varphi\in V).
\]
For fixed $(n,\nu)$, seek $c_n,p_n,a_n,o_n\in X_n$ and
$f_n,v_n\in H$. Suppressing $\nu$ in the notation, set
\[
b_n=\mathsf C(c_n),\qquad \omega_n=\mathsf O(o_n),\qquad
\mathcal F_n=\Ffun[\mathsf C(v_n)],
\]
and
\[
\mathcal T_n=M_*(c_n)
\bigl[\chi_a\nabla S_*(R_\nu a_n)+\chi_f\nabla R_\nu f_n\bigr].
\]
For every $\varphi\in X_n$, prescribe
\begin{equation}
\begin{aligned}
&((c_n)_t,\varphi)+D_c(\nabla c_n,\nabla\varphi)
=(\mathcal T_n,\nabla\varphi)\\
&\qquad+(r_cM_*(c_n)H_P(\omega_n)
-\eta b_n[1-\mathsf C(v_n)],\varphi),
\end{aligned}
\label{eq:global-galerkin-c}
\end{equation}
\begin{equation}
\begin{aligned}
&((p_n)_t,\varphi)+D_p(\nabla p_n,\nabla\varphi)
+\mu_p(p_n,\varphi)\\
&\qquad=(\sigma_pT(\theta_0+\theta_1b_n),\varphi),
\end{aligned}
\label{eq:global-galerkin-p}
\end{equation}
\begin{equation}
\begin{aligned}
&\tau_a((a_n)_t,\varphi)+D_a(\nabla a_n,\nabla\varphi)
+((\mu_a+\kappa_ab_n)a_n,\varphi)\\
&\qquad=(\sigma_aTH_H(\omega_n),\varphi),
\end{aligned}
\label{eq:global-galerkin-a}
\end{equation}
\begin{equation}
\begin{aligned}
&\tau_o((o_n)_t,\varphi)+D_o(\nabla o_n,\nabla\varphi)\\
&\qquad+(\gamma_TT h_{K_T}(\omega_n)
+\gamma_cb_n h_{K_c}(\omega_n),\varphi)\\
&\qquad=(\lambda\mathcal F_n(o_b-o_n),\varphi),
\end{aligned}
\label{eq:global-galerkin-o}
\end{equation}
where $h_K(s)=s/(K+s)$ for $s\geq0$. The unprojected ODEs are
\begin{align*}
(f_n)_t&=\rho_f(1-f_n)-\kappa_fp_n^+f_n,\\
(v_n)_t&=\eta b_n[1-\mathsf C(v_n)]
-\delta_v\{\mathsf C(v_n)-\Ffun[\mathsf C(v_n)]\}.
\end{align*}
The diffusive initial values are $P_n c_0,P_n p_0,P_n a_0,P_n o_0$;
the ODE initial values are $f_0,v_0$. Pressure is defined by
$\pi_n=\pi[\mathsf C(v_n)]$. Notice that the linear damping in
$p_n,a_n,o_n$ has not been truncated.

These equations form an ODE on
$\mathbb R^{4n}\times H\times H$ with a locally Lipschitz vector field
that is bounded on bounded sets. Indeed, norms on $X_n$ are equivalent,
$X_n\subset L^\infty$, and $R_\nu$ maps $X_n$ into itself.
For example,
\[
\|p^+f-\widetilde p^+\widetilde f\|_2
\leq\|p-\widetilde p\|_\infty\|f\|_2
+\|\widetilde p\|_\infty\|f-\widetilde f\|_2.
\]
The bounded, Lipschitz derivative $S_*'$ similarly controls differences
of $\nabla S_*(R_\nu a)$ on bounded subsets of $X_n$.
Theorem~\ref{thm:perfusion} and the $1$-Lipschitz property of
$\mathsf C$ give, for $0\leq b,\widetilde b\leq1$,
\[
\|\mathcal B_b(w)-\mathcal B_{\widetilde b}(\widetilde w)\|_2
\leq L_B\|w-\widetilde w\|_2+\eta\|b-\widetilde b\|_2,
\]
with $L_B=\eta+\delta_v(1+L_{\Ffun})$ and $\mathcal B$ as in
\eqref{eq:vessel-extension}. Products involving the pressure closure
and a spectral variable are locally Lipschitz by the same
$L^2$--$L^\infty$ estimate. The Banach-space Picard--Lindel\"of theorem
therefore gives a unique maximal approximate solution.

The ODE components already preserve their intervals. In fact,
\begin{equation*}
f_n(t)=f_0\exp\left[-\int_0^t(\rho_f+\kappa_fp_n^+(s))\,\mathrm ds\right]
+\rho_f\int_0^t
\exp\left[-\int_s^t(\rho_f+\kappa_fp_n^+(r))\,\mathrm dr\right]
\,\mathrm ds,
\end{equation*}
so $0\leq f_n\leq1$. For the vessel equation, the pointwise identity
\[
\Ffun[\mathsf C(w)](x)
=\mathsf C(w(x))G(q[\mathsf C(w)](x))
\]
shows that its right-hand side equals $\eta b_n\geq0$ on
$\{v_n<0\}$, and equals
$-\delta_v[1-G(q[\mathsf C(v_n)])]\leq0$ on $\{v_n>1\}$.
Testing with $-v_n^-$ and $(v_n-1)^+$ proves $0\leq v_n\leq1$.
Thus $\mathcal F_n=\Ffun[v_n]$ and $\pi_n=\pi[v_n]$.

\medskip\noindent
\textit{Step 2: Estimates before taking the Galerkin limit.}
All estimates in this step hold without assuming
$c_n,p_n,a_n,o_n\geq0$. Testing the protease equation by $p_n$ yields
\begin{equation}
\frac12\frac{\mathrm d}{\mathrm dt}\|p_n\|_2^2
+D_p\|\nabla p_n\|_2^2+\frac{\mu_p}{2}\|p_n\|_2^2\leq C.
\label{eq:p-energy}
\end{equation}
At this approximation level the explicit formula for $f_n$ and the
Sobolev chain rule imply $f_n(t)\in H^1$ and
\begin{equation*}
\partial_t\nabla f_n
=-(\rho_f+\kappa_fp_n^+)\nabla f_n
-\kappa_ff_n\nabla p_n^+.
\end{equation*}
Indeed, $p_n$ is bounded in space on compact subintervals of its
existence interval, and $\nabla p_n^+=\mathbf1_{\{p_n>0\}}\nabla p_n$.
Consequently,
\begin{equation}
\frac12\frac{\mathrm d}{\mathrm dt}\|\nabla f_n\|_2^2
+\frac{\rho_f}{2}\|\nabla f_n\|_2^2
\leq\frac{\kappa_f^2}{2\rho_f}\|\nabla p_n\|_2^2.
\label{eq:f-gradient-energy}
\end{equation}
This controls $f_n$ in $L^\infty(0,T_f;V)$.

Testing the VEGF and oxygen equations by $a_n$ and $o_n$, respectively,
gives
\begin{align}
\frac{\tau_a}{2}\frac{\mathrm d}{\mathrm dt}\|a_n\|_2^2
+D_a\|\nabla a_n\|_2^2+\frac{\mu_a}{2}\|a_n\|_2^2
&\leq C,
\label{eq:a-energy}\\
\frac{\tau_o}{2}\frac{\mathrm d}{\mathrm dt}\|o_n\|_2^2
+D_o\|\nabla o_n\|_2^2
+\frac{\lambda}{2}\int_\Omega\mathcal F_n o_n^2
&\leq\frac{\lambda o_b^2}{2}|\Omega|.
\label{eq:o-energy}
\end{align}
Here $b_n\geq0$, $\mathcal F_n\geq0$, and
$h_K(\mathsf O(s))s\geq0$ for every $s\in\mathbb R$.
The last fact justifies discarding the oxygen-consumption terms even
when $o_n$ is not nonnegative.

Since $|M_*|\leq1/4$, $|S_*'|\leq1$, and the endothelial reaction
in \eqref{eq:global-galerkin-c} is bounded, testing by $c_n$ and using
\eqref{eq:resolvent-estimates} gives
\begin{equation}
\begin{split}
\frac12\frac{\mathrm d}{\mathrm dt}\|c_n\|_2^2
+\frac{D_c}{2}\|\nabla c_n\|_2^2
\leq C\bigl(1+\|c_n\|_2^2
+\|\nabla a_n\|_2^2+\|\nabla f_n\|_2^2\bigr).
\end{split}
\label{eq:c-energy}
\end{equation}
Integrate first \eqref{eq:p-energy}--\eqref{eq:f-gradient-energy}
and \eqref{eq:a-energy}--\eqref{eq:o-energy}, and then apply
Gronwall's inequality to \eqref{eq:c-energy}. This yields
\begin{equation}
\begin{split}
&\sup_{0\leq t\leq T_f}
\left(\sum_{z_n=c_n,p_n,a_n,o_n}\|z_n(t)\|_2^2
+\|f_n(t)\|_V^2\right)\\
&\qquad+
\int_0^{T_f}\sum_{z_n=c_n,p_n,a_n,o_n}\|\nabla z_n\|_2^2
\,\mathrm dt\leq C_{T_f}.
\end{split}
\label{eq:global-energy}
\end{equation}
The ODEs also give
\[
\|(f_n)_t\|_2\leq\rho_f|\Omega|^{1/2}+\kappa_f\|p_n\|_2,
\qquad
\|(v_n)_t\|_2\leq(\eta+\delta_v)|\Omega|^{1/2}.
\]
These bounds prevent finite-time escape in
$\mathbb R^{4n}\times H\times H$. Since the vector field is bounded
on bounded sets, the approximate solution extends globally.

To estimate time derivatives in $V^*$, test the projected equations
against $P_n\varphi$, with $\varphi\in V$. The uniform bound
$\|P_n\varphi\|_V\leq\|\varphi\|_V$, the $L^2$ bounds on the taxis
fluxes, and the at-most-linear growth of the reactions imply
\begin{equation}
\sum_{z_n=c_n,p_n,a_n,o_n}
\|(z_n)_t\|_{L^2(0,T_f;V^*)}\leq C_{T_f}.
\label{eq:time-derivative-detail}
\end{equation}
All constants in \eqref{eq:global-energy}--\eqref{eq:time-derivative-detail}
are independent of both $n$ and $\nu$.

\medskip\noindent
\textit{Step 3: The limit $n\to\infty$ at fixed $\nu$.}
The compact embedding $V\Subset H$ and the continuous embedding
$H\hookrightarrow V^*$ allow the Aubin--Lions lemma to be applied
to each diffusive variable. After extraction,
\[
z_n\to z^\nu\quad\text{in }L^2(Q_T),\qquad
z_n\rightharpoonup z^\nu\quad\text{in }L^2(0,T_f;V),
\quad z=c,p,a,o,
\]
and a further subsequence converges almost everywhere.

The ODE variables converge strongly as well. Subtracting two ECM
equations, testing by $w=f_n-f_j$, and using $0\leq f_j\leq1$ gives
\[
\frac{\mathrm d}{\mathrm dt}\|w\|_2^2+\rho_f\|w\|_2^2
\leq\frac{\kappa_f^2}{\rho_f}\|p_n-p_j\|_2^2.
\]
The common initial value therefore implies
\[
\|f_n-f_j\|_{C([0,T_f];H)}^2
\leq C_{T_f}\|p_n-p_j\|_{L^2(Q_T)}^2.
\]
Similarly, Lemma~\ref{lem:vascular}, applied with drivers
$b_n=\mathsf C(c_n)$, gives
\[
\|v_n-v_j\|_{C([0,T_f];H)}^2
\leq C_{T_f}\|\mathsf C(c_n)-\mathsf C(c_j)\|_{L^2(Q_T)}^2.
\]
Thus $f_n\to f^\nu$ and $v_n\to v^\nu$ in $C([0,T_f];H)$;
the latter convergence does not use spatial regularity of $v_n$.
Theorem~\ref{thm:perfusion} also gives
\[
(\pi[v_n],q[v_n],\Ffun[v_n])
\longrightarrow
(\pi[v^\nu],q[v^\nu],\Ffun[v^\nu])
\quad\text{in }C([0,T_f];V\times H\times H).
\]

At fixed $\nu$, boundedness of $\nabla R_\nu:H\to L^2$ shows
that $\nabla R_\nu a_n\to\nabla R_\nu a^\nu$ and
$\nabla R_\nu f_n\to\nabla R_\nu f^\nu$ strongly in $L^2(Q_T)$.
The bounded multipliers $M_*(c_n)$ and $S_*'(R_\nu a_n)$ converge
almost everywhere, so the taxis fluxes converge strongly in $L^2$.
All reaction terms converge in $L^2(Q_T)$ by strong convergence and
boundedness of the clipped factors. For a typical term, split
\[
b_na_n-b^\nu a^\nu
=b_n(a_n-a^\nu)+(b_n-b^\nu)a^\nu;
\]
the first term converges in $L^2$, and dominated convergence handles
the second. The same argument applies to $\Ffun[v_n]o_n$ and
$p_n^+f_n$. Testing with $P_n\varphi$ and using
$P_n\varphi\to\varphi$ in $V$ now yields the regularized weak
system with all the truncations specified in Step~1.

For each diffusive limit,
$z^\nu\in L^2(0,T_f;V)$ and $(z^\nu)_t\in L^2(0,T_f;V^*)$.
The Hilbert-triple continuity theorem gives
$z^\nu\in C([0,T_f];H)$. Passing to the time-integrated weak
equations with test functions nonzero at $t=0$, and using
$P_nz_0\to z_0$ in $H$, identifies $z^\nu(0)=z_0$.
The ODE initial values follow from their convergence in $C([0,T_f];H)$.

\medskip\noindent
\textit{Step 4: Invariant intervals and comparison bounds.}
We may now use Sobolev truncations as test functions in the regularized
weak equations. In this step suppress the superscript $\nu$.
Testing the endothelial equation by $-c^-$ and $(c-1)^+$ gives
\[
\frac12\frac{\mathrm d}{\mathrm dt}\|c^-\|_2^2
+D_c\|\nabla c^-\|_2^2\leq0,
\qquad
\frac12\frac{\mathrm d}{\mathrm dt}\|(c-1)^+\|_2^2
+D_c\|\nabla(c-1)^+\|_2^2\leq0.
\]
Indeed, $M_*(c)=0$ on both excess sets, while the truncated reaction
is zero on $\{c<0\}$ and nonpositive on $\{c>1\}$.
Thus $0\leq c\leq1$. The nonnegative source and positive damping
in the protease and VEGF equations give $p,a\geq0$ by the tests
$-p^-$ and $-a^-$. Define
\[
P_*:=\max\left\{\|p_0\|_\infty,
\frac{\sigma_p(\theta_0+\theta_1)}{\mu_p}\right\},\qquad
A_*:=\max\left\{\|a_0\|_\infty,\frac{\sigma_a}{\mu_a}\right\}.
\]
Testing by $(p-P_*)^+$ and $(a-A_*)^+$ proves $p\leq P_*$ and
$a\leq A_*$. In particular, no bound on the $L^\infty$ norms of
the projected initial data was required.

For oxygen, the consumption terms vanish on $\{o<0\}$ because
$\mathsf O(o)=0$, and $\lambda\Ffun[v](o_b-o)\geq0$ there.
Testing by $-o^-$ proves $o\geq0$. On $\{o>o_b\}$, the
consumption terms are nonnegative and the perfusion term is
$-\lambda\Ffun[v](o-o_b)$. The test $(o-o_b)^+$ therefore proves
$o\leq o_b$. Together with the ODE bounds, this yields
\eqref{eq:invariant} and \eqref{eq:comparison-bounds}.
All reaction truncations are now inactive, $M_*(c)=M(c)$, and
$S_*(R_\nu a)=S(R_\nu a)$ by positivity of $R_\nu$.
Only the two resolvent regularizations in the taxis flux remain.

The estimates of Step~2 pass to this regularized solution by weak
lower semicontinuity. In particular, $f^\nu$ is bounded in $L^\infty(0,T_f;V)\cap W^{1,\infty}(0,T_f;H)$, $v^\nu$ is bounded in $W^{1,\infty}(0,T_f;H)$,
uniformly in $\nu$. The asserted $H^1$ bound for $f^\nu$ is thus
inherited from the justified approximate gradient estimate, without
formally differentiating an unconstructed weak solution.

\medskip\noindent
\textit{Step 5: Removal of the resolvent regularization.}
Choose $\nu_k\downarrow0$ and write $z_k=z^{\nu_k}$. Repeating
Aubin--Lions and the two ODE stability estimates gives, after extraction,
\begin{equation*}
\begin{aligned}
c_k,p_k,a_k,o_k&\to c,p,a,o &&\text{in }L^2(Q_T),\\
f_k,v_k&\to f,v &&\text{in }C([0,T_f];H),
\end{aligned}
\end{equation*}
with weak convergence of the diffusive variables in $L^2(0,T_f;V)$
and weak-star convergence of $f_k$ in $L^\infty(0,T_f;V)$.
Their time-derivative bounds also pass to the limit.
Theorem~\ref{thm:perfusion} gives strong convergence of
$\pi[v_k],q[v_k],\Ffun[v_k]$ in $C([0,T_f];V\times H\times H)$.
The interval and comparison bounds are retained.

For $z_k=a_k$ or $f_k$, \eqref{eq:resolvent-estimates} gives
\[
\|R_{\nu_k}z_k-z_k\|_{L^2(Q_T)}
\leq\sqrt{\nu_k}\,\|\nabla z_k\|_{L^2(Q_T)}\longrightarrow0.
\]
Consequently, $R_{\nu_k}a_k\to a,\qquad R_{\nu_k}f_k\to f$ in $L^2(Q_T)$,
and their gradients converge weakly in $L^2(Q_T)$ to
$\nabla a$ and $\nabla f$, respectively. This gradient identification
follows by testing distributional derivatives and using the strong
$L^2$ convergence.

Set
\[
B_k=M(c_k)S'(R_{\nu_k}a_k),\qquad B=M(c)S'(a).
\]
After extraction, $B_k\to B$ almost everywhere and $|B_k|\leq1/4$.
For every fixed $\Psi\in L^2(Q_T;\mathbb R^2)$, dominated
convergence gives $B_k\Psi\to B\Psi$ in $L^2$. Hence
\[
\begin{aligned}
\int_{Q_T}B_k\nabla R_{\nu_k}a_k\cdot\Psi =\int_{Q_T}\nabla R_{\nu_k}a_k\cdot(B_k\Psi)
\longrightarrow\int_{Q_T}\nabla a\cdot(B\Psi).
\end{aligned}
\]
Thus
\[
M(c_k)\nabla S(R_{\nu_k}a_k)
\rightharpoonup M(c)S'(a)\nabla a
\quad\text{in }L^2(Q_T;\mathbb R^2).
\]
The same argument with multiplier $M(c_k)$ proves
\[
M(c_k)\nabla R_{\nu_k}f_k
\rightharpoonup M(c)\nabla f
\quad\text{in }L^2(Q_T;\mathbb R^2).
\]
The Sobolev chain rule identifies $S'(a)\nabla a=\nabla S(a)$.
All reaction terms converge in $L^2(Q_T)$, since the state variables
are bounded and converge strongly, and the functional density
converges strongly by Theorem~\ref{thm:perfusion}. In particular,
\[
\Ffun[v_k](o_b-o_k)\to\Ffun[v](o_b-o),\qquad
v_k-\Ffun[v_k]\to v-\Ffun[v]
\quad\text{in }L^2(Q_T).
\]
Passing to the weak formulations now gives
\eqref{eq:weak-c-compact}--\eqref{eq:weak-o-compact}, while
$\pi=\pi[v]$ satisfies \eqref{eq:weak-pressure}.

Finally, $c,p,a,o\in L^2(0,T_f;V)$ with time derivatives in
$L^2(0,T_f;V^*)$, so they belong to $C([0,T_f];H)$.
The time-integrated weak formulations identify their initial values
as $c_0,p_0,a_0,o_0$; the strong $C([0,T_f];H)$ convergence identifies
the ODE initial values. The uniform estimates give exactly
\eqref{eq:weak-reg-diffusive}--\eqref{eq:weak-reg-pressure}, and
Theorem~\ref{thm:perfusion} gives
$q[v]\in L^\infty(0,T_f;L^2(\Omega))$.

The approximate trajectories were constructed globally. Performing
the preceding subsequence extractions successively on
$[0,1],[0,2],\ldots$, first at fixed $\nu$ and then along
$\nu_k\downarrow0$, gives a single global weak solution by a diagonal
argument. No uniqueness of the coupled weak system is needed.
\end{proof}

\section{Fast-signal limit}
\label{app:fast-details}
We provide the full proof of Theorem~\ref{thm:fast}. 

\begin{proof}
Fix $T_f>0$ and write
\[
H=L^2(\Omega),\qquad V=H^1(\Omega),\qquad
\mathcal X=L^2(0,T_f;V),\qquad
\mathcal X^*=L^2(0,T_f;V^*).
\]
Here $(\cdot,\cdot)$ denotes the $H$ inner product, and
$Q_T=\Omega\times(0,T_f)$. After discarding finitely many terms, we
may assume $0<\tau_o^n,\tau_a^n\leq1$. Constants may depend on $T_f$,
the common initial data, and the fixed parameters, but not on $n$.

\medskip\noindent
\textit{Step 1: Estimates uniform in both relaxation times.}
Set
\[
P_*:=\max\left\{\|p_0\|_\infty,
\frac{\sigma_p(\theta_0+\theta_1)}{\mu_p}\right\},\qquad
A_*:=\max\left\{\|a_0\|_\infty,\frac{\sigma_a}{\mu_a}\right\}.
\]
The invariant-region and comparison estimates give
\begin{equation}
0\leq c_n,f_n,v_n\leq1,\qquad
0\leq o_n\leq o_b,\qquad
0\leq p_n\leq P_*,\qquad 0\leq a_n\leq A_*.
\label{eq:fast-interval-bounds}
\end{equation}
These bounds are independent of both relaxation times.

For clarity, the signal energy estimates are obtained without dividing
by $\tau_a^n$ or $\tau_o^n$. Testing their equations by $a_n,o_n$
and integrating in time yields
\begin{align*}
D_a\int_0^{T_f}\|\nabla a_n\|_2^2\,\mathrm dt
+\frac{\mu_a}{2}\int_0^{T_f}\|a_n\|_2^2\,\mathrm dt
&\leq\frac{\tau_a^n}{2}\|a_0\|_2^2
+\frac{\sigma_a^2|\Omega|T_f}{2\mu_a},\\
D_o\int_0^{T_f}\|\nabla o_n\|_2^2\,\mathrm dt
&\leq\frac{\tau_o^n}{2}\|o_0\|_2^2
+\frac{\lambda o_b^2|\Omega|T_f}{2}.
\end{align*}
The oxygen estimate uses $0\leq\Ffun[v_n]\leq1$ and retains the
nonnegative absorption terms before applying Young's inequality.

The protease energy estimate gives a uniform bound for
$p_n$ in $L^\infty(0,T_f;H)\cap\mathcal X$. The explicit ECM ODE,
or its spatial difference quotients, gives
\[
\frac12\frac{\mathrm d}{\mathrm dt}\|\nabla f_n\|_2^2
+\frac{\rho_f}{2}\|\nabla f_n\|_2^2
\leq\frac{\kappa_f^2}{2\rho_f}\|\nabla p_n\|_2^2.
\]
Thus $f_n$ is uniformly bounded in $L^\infty(0,T_f;V)$.
Testing the endothelial equation by $c_n$, using $M(c_n)\leq1/4$
and $|S'(a_n)|\leq1$, gives
\[
\frac12\frac{\mathrm d}{\mathrm dt}\|c_n\|_2^2
+\frac{D_c}{2}\|\nabla c_n\|_2^2
\leq C\bigl(1+\|\nabla a_n\|_2^2+\|\nabla f_n\|_2^2\bigr).
\]
Together with \eqref{eq:fast-interval-bounds}, these inequalities imply
\begin{equation}
\|o_n\|_{\mathcal X}+\|a_n\|_{\mathcal X}
+\|c_n\|_{\mathcal X}+\|p_n\|_{\mathcal X}
+\|f_n\|_{L^\infty(0,T_f;V)}\leq C_{T_f}.
\label{eq:fast-uniform}
\end{equation}
The slow equations and the ODE bounds further imply
\begin{equation*}
\begin{aligned}
\|(c_n)_t\|_{\mathcal X^*}+\|(p_n)_t\|_{\mathcal X^*}
&\leq C_{T_f},\\
\|(f_n)_t\|_{L^\infty(0,T_f;H)}
+\|(v_n)_t\|_{L^\infty(0,T_f;H)}&\leq C_{T_f}.
\end{aligned}
\end{equation*}
No uniform estimate for the unscaled derivatives $(o_n)_t,(a_n)_t$
is needed.

\medskip\noindent
\textit{Step 2: Compactness of the slow variables and perfusion.}
Since $V\Subset H\hookrightarrow V^*$, Aubin--Lions gives, after
extraction,
\[
c_n,p_n\to c,p\quad\text{in }L^2(Q_T),\qquad
c_n,p_n\rightharpoonup c,p\quad\text{in }\mathcal X.
\]
Subtracting two ECM equations and testing by $f_n-f_j$ yields
\[
\frac{\mathrm d}{\mathrm dt}\|f_n-f_j\|_2^2
+\rho_f\|f_n-f_j\|_2^2
\leq\frac{\kappa_f^2}{\rho_f}\|p_n-p_j\|_2^2.
\]
The common initial value therefore gives
\[
\|f_n-f_j\|_{C([0,T_f];H)}^2
\leq C_{T_f}\|p_n-p_j\|_{L^2(Q_T)}^2.
\]
Likewise, Lemma~\ref{lem:vascular} with drivers $c_n,c_j$ gives
\[
\|v_n-v_j\|_{C([0,T_f];H)}^2
\leq C_{T_f}\|c_n-c_j\|_{L^2(Q_T)}^2.
\]
Hence
\[
f_n\to f,\qquad v_n\to v
\quad\text{in }C([0,T_f];H).
\]
After extraction, $f_n\rightharpoonup^*f$ in $L^\infty(0,T_f;V)$,
and $(f_n)_t,(v_n)_t$ converge weak-star to $f_t,v_t$ in
$L^\infty(0,T_f;H)$.
Theorem~\ref{thm:perfusion} implies
\begin{equation}
\begin{aligned}
\pi[v_n]&\to\pi[v]=:\pi
&&\text{in }C([0,T_f];V),\\
q[v_n]&\to q[v]
&&\text{in }C([0,T_f];H),\\
F_n:=\Ffun[v_n]&\to F:=\Ffun[v]
&&\text{in }C([0,T_f];H).
\end{aligned}
\label{eq:fast-closure-convergence}
\end{equation}
In particular, \eqref{eq:fast-perfusion} holds.

Further extraction gives $o_n\rightharpoonup o$, $a_n\rightharpoonup a$ in $\mathcal X$.
All bounded components converge weak-star in $L^\infty(Q_T)$.
In particular,
\[
0\leq o\leq o_b,\qquad 0\leq a\leq A_*.
\]
We also arrange that $c_n\to c$ and $F_n\to F$ almost everywhere
in $Q_T$; both sequences take values in $[0,1]$.

\medskip\noindent
\textit{Step 3: Scaled residuals and the oxygen operator.}
Define
\[
r_n^o:=\tau_o^n(o_n)_t,\qquad
r_n^a:=\tau_a^n(a_n)_t.
\]
The signal equations, \eqref{eq:fast-uniform}, and the interval bounds
give $\|r_n^o\|_{\mathcal X^*}+\|r_n^a\|_{\mathcal X^*}\leq C_{T_f}$.
For every $\zeta\in C_c^\infty(0,T_f;V)$,
\[
\int_0^{T_f}\langle r_n^o,\zeta\rangle\,\mathrm dt
=-\tau_o^n\int_0^{T_f}(o_n,\zeta_t)\,\mathrm dt\longrightarrow0,
\]
and the analogous identity holds for $r_n^a$. Density of these test
functions in $\mathcal X$ and the uniform dual bounds imply
\[
r_n^o\rightharpoonup0,\qquad r_n^a\rightharpoonup0
\quad\text{in }\mathcal X^*,
\]
which proves \eqref{eq:fast-residuals}.

For $K>0$, write $d_K=K/(K+o_b)^2$ and use the extension
\[
\widehat h_K(s)=
\begin{cases}
s/K,&s<0,\\
s/(K+s),&0\leq s\leq o_b,\\
o_b/(K+o_b)+d_K(s-o_b),&s>o_b.
\end{cases}
\]
It is globally Lipschitz, vanishes at zero, and satisfies
\[
(\widehat h_K(r)-\widehat h_K(s))(r-s)
\geq d_K|r-s|^2\qquad(r,s\in\mathbb R).
\]
Define, for almost every $t$,
\begin{equation*}
\begin{aligned}
\langle\mathcal L_n^oz,\varphi\rangle
={}&D_o(\nabla z,\nabla\varphi)
+\gamma_T(T\widehat h_{K_T}(z),\varphi)\\
&+\gamma_c(c_n\widehat h_{K_c}(z),\varphi)
+\lambda(F_nz,\varphi),
\end{aligned}
\end{equation*}
and let $\mathcal L^o$ be the same operator with $c,F$ in place of
$c_n,F_n$. We use the same notation for the induced maps
$\mathcal X\to\mathcal X^*$. They are globally Lipschitz, with a
constant independent of $n$.

The weighted Poincar\'e inequality of
Lemma~\ref{lem:weighted-poincare} gives, for $z,w\in V$,
\begin{equation}
\begin{aligned}
\langle\mathcal L_n^oz-\mathcal L_n^ow,z-w\rangle
&\geq D_o\|\nabla(z-w)\|_2^2
+\gamma_Td_{K_T}\int_\Omega T|z-w|^2\\
&\geq\alpha_o\|z-w\|_V^2,
\end{aligned}
\label{eq:fast-oxygen-coercivity}
\end{equation}
where one may take
\[
\alpha_o=\frac{\min\{D_o,\gamma_Td_{K_T}\}}
{1+C_{\mathrm{wp}}}>0.
\]
This constant is independent of $n$ and requires no positive lower
bound on $F_n$.

For every fixed $z\in\mathcal X$, dominated convergence gives
\[
(c_n-c)\widehat h_{K_c}(z)\to0,\qquad
(F_n-F)z\to0
\quad\text{in }L^2(Q_T).
\]
Indeed, the coefficient differences are uniformly bounded and
converge almost everywhere, while $z,\widehat h_{K_c}(z)\in L^2(Q_T)$.
Consequently,
\begin{equation}
\mathcal L_n^oz\to\mathcal L^oz
\quad\text{in }\mathcal X^*
\qquad\text{for every fixed }z\in\mathcal X.
\label{eq:operator-coefficient-convergence}
\end{equation}

\medskip\noindent
\textit{Step 4: Oxygen identification and strong convergence.}
Set $g_n^o=\lambda o_bF_n$ and $g^o=\lambda o_bF$. Then
\[
r_n^o+\mathcal L_n^oo_n=g_n^o,
\qquad
g_n^o\to g^o\quad\text{in }L^2(Q_T),
\]
and therefore $\mathcal L_n^oo_n\rightharpoonup g^o$ in
$\mathcal X^*$.
For each $n$, the Hilbert-triple chain rule permits testing the
dynamic equation by $o_n$. It yields
\begin{equation*}
\int_0^{T_f}\langle\mathcal L_n^oo_n,o_n\rangle\,\mathrm dt
=\int_0^{T_f}(g_n^o,o_n)\,\mathrm dt -\frac{\tau_o^n}{2}
\bigl(\|o_n(T_f)\|_2^2-\|o_0\|_2^2\bigr).
\end{equation*}
The endpoint values are well defined in $H$, and their norms are
bounded by $o_b|\Omega|^{1/2}$. Hence the last term tends to zero.
Strong convergence of $g_n^o$ and weak convergence of $o_n$ imply
\begin{equation}
\int_0^{T_f}\langle\mathcal L_n^oo_n,o_n\rangle\,\mathrm dt
\longrightarrow\int_0^{T_f}(g^o,o)\,\mathrm dt.
\label{eq:oxygen-pairing-limit}
\end{equation}

For every fixed $z\in\mathcal X$, integrated monotonicity gives
\[
\begin{aligned}
0\leq{}&\int_0^{T_f}
\langle\mathcal L_n^oo_n-\mathcal L_n^oz,o_n-z\rangle\,\mathrm dt\\
={}&\int_0^{T_f}\langle\mathcal L_n^oo_n,o_n\rangle\,\mathrm dt
-\int_0^{T_f}\langle\mathcal L_n^oo_n,z\rangle\,\mathrm dt\\
&-\int_0^{T_f}\langle\mathcal L_n^oz,o_n\rangle\,\mathrm dt
+\int_0^{T_f}\langle\mathcal L_n^oz,z\rangle\,\mathrm dt.
\end{aligned}
\]
Using \eqref{eq:oxygen-pairing-limit} for the first term,
$\mathcal L_n^oo_n\rightharpoonup g^o$ for the second, and
\eqref{eq:operator-coefficient-convergence} for the last two, we obtain
\begin{equation*}
\int_0^{T_f}\langle g^o-\mathcal L^oz,o-z\rangle\,\mathrm dt
\geq0\qquad(z\in\mathcal X).
\end{equation*}
For an arbitrary time-dependent $\varphi\in\mathcal X$, choose
$z=o-s\varphi$ and $z=o+s\varphi$, with $s>0$. Divide the resulting
inequalities by $s$ and let $s\downarrow0$. The Lipschitz continuity
of $\mathcal L^o:\mathcal X\to\mathcal X^*$ gives
\[
\int_0^{T_f}\langle\mathcal L^oo-g^o,\varphi\rangle\,\mathrm dt=0
\qquad(\varphi\in\mathcal X).
\]
Thus $\mathcal L^oo=g^o$ in $\mathcal X^*$, equivalently in $V^*$
for almost every time. Since $0\leq o\leq o_b$, the extension
agrees with the original oxygen nonlinearity.
Proposition~\ref{prop:resolvents}
therefore identifies $o=\mathcal O(c,v)$ almost everywhere.

To strengthen the convergence, apply
\eqref{eq:fast-oxygen-coercivity} with $z=o_n,w=o$ and integrate:
\begin{equation*}
\begin{aligned}
\alpha_o\|o_n-o\|_{\mathcal X}^2
\leq{}&\int_0^{T_f}\langle\mathcal L_n^oo_n,o_n\rangle\,\mathrm dt
-\int_0^{T_f}\langle\mathcal L_n^oo_n,o\rangle\,\mathrm dt\\
&-\int_0^{T_f}\langle\mathcal L_n^oo,o_n-o\rangle\,\mathrm dt
\longrightarrow0.
\end{aligned}
\end{equation*}
Indeed, the first two terms converge to the same value
$\int_0^{T_f}(g^o,o)\,\mathrm dt$, and the last tends to zero because
$\mathcal L_n^oo\to\mathcal L^oo$ strongly in $\mathcal X^*$ and
$o_n-o\rightharpoonup0$ in $\mathcal X$. Hence
\begin{equation}
o_n\to o\quad\text{strongly in }\mathcal X.
\label{eq:oxygen-strong}
\end{equation}

\medskip\noindent
\textit{Step 5: VEGF identification and strong convergence.}
Define
\begin{equation*}
\langle\mathcal L_n^az,\xi\rangle
=D_a(\nabla z,\nabla\xi)+((\mu_a+\kappa_ac_n)z,\xi),
\end{equation*}
and let $\mathcal L^a$ have coefficient $c$ in place of $c_n$.
Both act from $\mathcal X$ to $\mathcal X^*$, and
\[
\langle\mathcal L_n^az,z\rangle
\geq\alpha_a\|z\|_V^2,\qquad
\alpha_a:=\min\{D_a,\mu_a\}>0.
\]
Set
\[
g_n^a=\sigma_aTH_H(o_n),\qquad g^a=\sigma_aTH_H(o).
\]
By \eqref{eq:oxygen-strong} and the Lipschitz continuity of $H_H$,
$g_n^a\to g^a$ strongly in $L^2(Q_T)$. The equation reads $r_n^a+\mathcal L_n^aa_n=g_n^a$.
For every fixed $z\in\mathcal X$, dominated convergence gives
$\mathcal L_n^az\to\mathcal L^az$ strongly in $\mathcal X^*$.
Moreover, for every $\xi\in\mathcal X$,
\[
(\mu_a+\kappa_ac_n)\xi\to(\mu_a+\kappa_ac)\xi
\quad\text{in }L^2(Q_T).
\]
Pairing with $a_n\rightharpoonup a$ in $L^2(Q_T)$, and using weak
convergence of the gradients, therefore shows $\mathcal L_n^aa_n\rightharpoonup\mathcal L^aa
\quad\text{in }\mathcal X^*$.
Since $r_n^a\rightharpoonup0$, the limiting equation is
$\mathcal L^aa=g^a$. Thus $a=\mathcal A(c,o)$ almost everywhere,
and Proposition~\ref{prop:resolvents} gives the sharper elliptic bound
$0\leq a\leq\sigma_a/\mu_a$.

The energy identity for $a_n$ is
\begin{equation*}
\int_0^{T_f}\langle\mathcal L_n^aa_n,a_n\rangle\,\mathrm dt
=\int_0^{T_f}(g_n^a,a_n)\,\mathrm dt-\frac{\tau_a^n}{2}
\bigl(\|a_n(T_f)\|_2^2-\|a_0\|_2^2\bigr).
\end{equation*}
By \eqref{eq:fast-interval-bounds}, the endpoint norms are bounded
by $A_*|\Omega|^{1/2}$, so the boundary term vanishes. Consequently,
\[
\int_0^{T_f}\langle\mathcal L_n^aa_n,a_n\rangle\,\mathrm dt
\longrightarrow\int_0^{T_f}(g^a,a)\,\mathrm dt.
\]
Coercivity now gives
\[
\alpha_a\|a_n-a\|_{\mathcal X}^2
\leq\int_0^{T_f}\langle\mathcal L_n^aa_n,a_n\rangle\,\mathrm dt
-\int_0^{T_f}\langle\mathcal L_n^aa_n,a\rangle\,\mathrm dt-\int_0^{T_f}\langle\mathcal L_n^aa,a_n-a\rangle\,\mathrm dt
\longrightarrow0.
\]
The first two terms have the same limit, while the last tends to zero
by strong convergence of $\mathcal L_n^aa$ in $\mathcal X^*$ and
weak convergence of $a_n-a$ in $\mathcal X$. Hence $a_n\to a$ strongly in $\mathcal X$.
This proves \eqref{eq:fast-ao}.

\medskip\noindent
\textit{Step 6: The limiting slow equations and initial traces.}
After further extraction, the strongly convergent fields converge
almost everywhere in $Q_T$. Put
$B_n=M(c_n)S'(a_n)$ and $B=M(c)S'(a)$. Then
$B_n\to B$ almost everywhere and $|B_n|\leq1/4$. Thus
\[
\begin{aligned}
\|B_n\nabla a_n-B\nabla a\|_{L^2(Q_T)}
&\leq\tfrac14\|\nabla a_n-\nabla a\|_{L^2(Q_T)}
+\|(B_n-B)\nabla a\|_{L^2(Q_T)}\longrightarrow0,
\end{aligned}
\]
where the second term is handled by dominated convergence.
For haptotaxis, $\nabla f_n\rightharpoonup\nabla f$ in $L^2(Q_T)$,
and for every $\Psi\in L^2(Q_T;\mathbb R^2)$,
$M(c_n)\Psi\to M(c)\Psi$ strongly in $L^2$. Hence $M(c_n)\nabla f_n\rightharpoonup M(c)\nabla f$ in $L^2(Q_T;\mathbb R^2)$.
All slow reaction terms converge strongly in $L^2(Q_T)$ by the
strong convergence of the bounded state variables and
\eqref{eq:fast-closure-convergence}. For the vessel loss, one can use
the identity $[1-G(q[v_n])]v_n=v_n-F_n$.
Passing to the weak formulations therefore yields
\eqref{eq:weak-c-compact}, \eqref{eq:weak-p-compact},
\eqref{eq:weak-f-compact}, and \eqref{eq:weak-v-compact}.
The signal identifications give
\eqref{eq:elliptic-o-compact}--\eqref{eq:elliptic-a-compact}
in $V^*$ for almost every time, and $\pi=\pi[v]$ satisfies
\eqref{eq:weak-pressure}.

The uniform bounds on $c_n,p_n$ and their time derivatives imply
$c,p\in\mathcal X$ with $c_t,p_t\in\mathcal X^*$, and hence
$c,p\in C([0,T_f];H)$. Passing to the time-integrated slow equations
with test functions nonzero at $t=0$ identifies
$c(0)=c_0$ and $p(0)=p_0$. The strong convergence in
$C([0,T_f];H)$ gives $f(0)=f_0$ and $v(0)=v_0$.
The ODE derivative bounds pass weak-star in $L^\infty(0,T_f;H)$,
and all remaining regularity and invariant bounds are inherited
from the estimates above. No initial condition is imposed on the
limiting elliptic signals.

The strong convergences also give the asserted almost-everywhere
convergences after extraction; in particular, the fast gradients
converge almost everywhere along a further subsequence. A diagonal
extraction on $[0,1],[0,2],\ldots$ gives a single subsequence and a
global quasi-static weak solution. No relation between the rates
at which $\tau_o^n$ and $\tau_a^n$ vanish, no compatibility of the
initial signals with the elliptic equations, and no uniqueness of
the full limiting system have been used.
\end{proof}

\end{appendices}

\end{document}